\documentclass[11pt]{article}
\usepackage[margin=1in]{geometry}
\usepackage{amsmath,amssymb,amsthm,mathtools}
\usepackage{graphicx}
\usepackage{booktabs}
\usepackage{longtable}
\usepackage{enumitem}
\usepackage{needspace}
\usepackage{microtype}
\usepackage[section]{placeins}
\usepackage[hidelinks]{hyperref}
\usepackage[capitalise,nameinlink]{cleveref}

\newtheorem{theorem}{Theorem}[section]
\newtheorem{proposition}[theorem]{Proposition}
\newtheorem{lemma}[theorem]{Lemma}
\newtheorem{corollary}[theorem]{Corollary}
\newtheorem{conjecture}[theorem]{Conjecture}
\theoremstyle{definition}
\newtheorem{definition}[theorem]{Definition}

\newtheorem{problem}[theorem]{Problem}
\theoremstyle{remark}
\newtheorem{remark}[theorem]{Remark}
\numberwithin{equation}{section}
\numberwithin{table}{section}
\numberwithin{figure}{section}

\newcommand{\N}{\mathbb{N}}

\newcommand{\Zt}{\mathbb{Z}_2}
\newcommand{\Zu}{\mathbb{Z}_2^{\times}}
\newcommand{\cS}{\mathcal{S}}
\newcommand{\E}{\mathbb{E}}
\renewcommand{\P}{\mathbb{P}}

\newcommand{\pair}[2]{\langle #1,\,#2\rangle}
\DeclareMathOperator{\Syr}{Syr}
\newcommand{\even}{2\N}
\newcommand{\odd}{2\N-1}
\newcommand{\So}{S_{\mathrm o}}
\newcommand{\Se}{S_{\mathrm e}}

\title{Up or down by half: a family of Collatz-type maps}
\author{Oleg Vasilyev\\
[2pt] {\small California, USA}\\ 
{\small\texttt{vasilyev.oleg@gmail.com}}}
\date{}

\begin{document}

\maketitle

\begin{abstract}
Each step of the Collatz map changes a positive integer by about half of its value: odd numbers go up and even numbers go down. We keep the size of the steps and vary their direction. For a set $S$ of positive integers and $x=2^m k$ with $k$ odd, the map $T_S$ sends $x$ to $(3x+1)/2$ if $m=0$, to $3x/2$ if $m\in S$, and to $k$ otherwise; the map $T^-_S$ differs only at odd $x$, which it sends down to $(x+1)/2$. The Collatz map is $T_\emptyset$ (with runs of halvings merged) and, after a shift by one, $T^-_{\mathbb N}$. Independence of the 2-adic valuations along odd-to-odd steps gives random-walk models with drifts $(1+\rho(S))\log_2 3-2$ for $T_S$ and $\rho(S)\log_2 3-2$ for $T^-_S$, where $\rho(S)$ is an explicit series. Thus every $T^-_S$ is contracting, while $T_S$ is contracting exactly when $\rho(S)<\log_3(4/3)\approx 0.2619$, and expanding otherwise. We prove the corresponding density results, classify the intervals $S$, and give a simple criterion under which every orbit of $T^-_S$ reaches $1$ or the cycle $(2,3)$; beyond it, residue classes alone cannot prove that all orbits descend. For seven contracting maps $T_S$, three with $\rho=1/4$ just below the threshold, every $n<2^{34}$ ($n<1.668\cdot 10^{10}$ for $S=\{1\}$) reaches one of at most five cycles, and the stopping-time statistics agree with the model to three to five digits. We then introduce branching maps, in which any member of a family may be applied at each step, and note a transfer principle: if a family is bounded, so is every member; e.g., boundedness of $\{\emptyset,\{3\}\}$ would imply that no Collatz orbit diverges. For $\{\emptyset,\{1\}\}$, a family of two contracting maps, we prove that a rule depending only on $x\bmod 32$ has drift $\tfrac{10}{17}\log_2\tfrac98>0$ for almost every 2-adic starting point, deduce that for almost all $n$ the reachable set contains numbers larger than $n^{1.05}$, and show that residues modulo $2^L$, $L\le 10$, cannot prove that this family is unbounded. Finally, we let a valuation decide the direction at odd numbers as well. In the resulting unified maps the odd $x\equiv 1\pmod 4$ decide the behavior: if they go down, every block of steps is dominated by the corresponding block of the Collatz map, and a census finds only the cycles $(1)$, $(2,3)$ and $(3,5)$ with smallest element below $10^9$; if they go up, the maps behave like the maps $T_S$.
\end{abstract}

\medskip
\noindent\textbf{Keywords:} Collatz problem; $3x+1$ problem; generalized Collatz maps; random walks;
branching random walks; 2-adic integers.

\noindent\textbf{2020 MSC:} 11B83; 11Y55; 60G50; 60J80.

\section{Introduction}\label{sec:intro}

The Collatz map
\[
T(x)=\begin{cases}(3x+1)/2, & x \text{ odd},\\ x/2, & x\text{ even},\end{cases}
\]
changes a positive integer by about half of its value at every step: an odd $x$ goes up by
$(x+1)/2$, an even $x$ goes down by $x/2$. The Collatz conjecture states that every
orbit eventually reaches the cycle $(1,2)$; see Lagarias's surveys and annotated bibliographies
\cite{Lag85,Lag10,LagOverview,LagBib1,LagBib2} and the monograph of Wirsching \cite{Wir98}. In this
paper we keep the size of the steps and vary their direction: at some numbers we go up where the
Collatz map goes down, or down where it goes up. This gives a family of maps that contains the
Collatz map, maps whose orbits are expected to diverge, and maps for which the analogue of the
conjecture is easy to prove. One can then ask which members are easier than the Collatz map, and
whether statements can be transferred from one member to another.

We let the direction depend on the 2-adic valuation of the current number. Write every positive
integer uniquely as $x=2^m k$ with $k$ odd, and let $\N=\{1,2,3,\dots\}$.

\begin{definition}\label{def:TS}
For $S\subseteq\N$ let $T_S,T^-_S\colon\N\to\N$ be
\begin{equation}\label{eq:TS}
T_S(x)=\begin{cases}
(3x+1)/2, & m=0,\\
3x/2, & m\in S,\\
k, & m\ge 1,\ m\notin S,
\end{cases}
\qquad
T^-_S(x)=\begin{cases}
(x+1)/2, & m=0,\\
3x/2, & m\in S,\\
k, & m\ge 1,\ m\notin S,
\end{cases}
\end{equation}
where $x=2^m k$ with $k$ odd.
\end{definition}

Thus $S$ is the set of exponents $m\ge 1$ at which the map goes up, and odd numbers go up under
$T_S$ and down under $T^-_S$; equivalently, the two families together are indexed by the set of
all exponents $m\ge 0$ at which the map goes up. A downward step at an even number removes all
factors of $2$ at once. For the Collatz map this only merges a run of halvings into one step:
$T_\emptyset$ sends an even number directly to its odd part, so its orbits are the Collatz orbits
with every run of halvings compressed, and $T_\emptyset$ has a divergent or a nontrivial periodic
orbit exactly when $T$ does. (\Cref{rem:single} discusses the variant in which a downward step
halves only once.)

For example, for $S=\{1\}$ the number $26=2\cdot 13$ has $m=1\in S$ and goes up to $39$, while
$680=2^3\cdot 85$ has $m=3\notin S$ and goes down to $85$; under the Collatz map $T_\emptyset$ both
would go down to their odd parts $13$ and $85$.

An orbit of a map $f\colon\N\to\N$ is either eventually periodic, in which case it ends in a
\emph{cycle} $\{x,f(x),\dots,f^{p-1}(x)\}$ with $f^p(x)=x$, or it consists of pairwise distinct
numbers and then tends to infinity; in the second case we say that the orbit \emph{diverges}. The
analogue of the Collatz conjecture for $T_S$ is the statement that no orbit of $T_S$ diverges.
Unlike the Collatz map, most maps $T_S$ have several cycles (see \cref{tab:cycles,tab:census}); we
write a cycle as a tuple that starts with its smallest element.

The Collatz map appears in the family a second time: $T^-_\N(x)=T(x-1)+1$ (\cref{prop:shift}).
So the family of maps $T_S$ starts at the Collatz map and its downward pull weakens as $S$ grows,
up to $T_\N(x)=\lceil 3x/2\rceil>x$, whose orbits all diverge; the family of maps $T^-_S$ ends at the
Collatz map and its pull strengthens as $S$ shrinks, down to $T^-_\emptyset$, whose orbits all
decrease to $1$. For example,
\begin{align*}
T_{\{1\}}&:\quad 7,\,11,\,17,\,26,\,39,\,59,\,89,\,134,\,201,\,302,\,453,\,680,\,85,\,128,\,1,\,2,\,3,\,5,\,8,\,1,\dots\\
T_{\{2,4,6,\dots\}}&:\quad 7,\,11,\,17,\,26,\,13,\,20,\,30,\,15,\,23,\,35,\,53,\,80,\,120,\,15,\dots\\
T^-_{\{1,2,3\}}&:\quad 15,\,8,\,12,\,18,\,27,\,14,\,21,\,11,\,6,\,9,\,5,\,3,\,2,\,3,\dots\\
T^-_{\N}&:\quad 28,\,42,\,63,\,32,\,48,\,72,\,108,\,162,\,243,\,122,\,183,\,92,\dots
\end{align*}
The first orbit ends in the cycle $(1,2,3,5,8)$, the second in $(15,23,35,53,80,120)$; the third and fourth end in $(2,3)$ (the fourth is the Collatz orbit of $27$ shifted by one; it takes 71 steps to first repeat). The questions that motivate this paper are: for which $S$ is the analogue of the Collatz conjecture plausible, what do the orbits of the maps near the boundary look like, and can a statement be proved more easily for one map than for another, or be transferred from one map to another?

\subsection*{The heuristic}
The following heuristic, standard for the Collatz map \cite{Cra78,LW92,KL10}, indicates which maps
$T_S$ should have only eventually periodic orbits. Consider only the odd numbers
$x_0,x_1,x_2,\dots$ of an orbit of $T_S$. The step from $x_i$ to $x_{i+1}$ begins with
$x\mapsto 3x+1$, which produces $M_i+1=v_2(3x_i+1)\ge 1$ factors $2$; the map then trades some of
them, say $r_i$, for factors $3$ by steps $x\mapsto 3x/2$, and removes the remaining ones. Thus
\[
x_{i+1}\approx\frac{3^{1+r_i}}{2^{M_i+1}}\,x_i,\qquad\text{that is,}\qquad
\log_2 x_{i+1}-\log_2 x_i\approx(1+r_i)\log_2 3-(M_i+1).
\]
For a random odd integer $x$ the valuation $v_2(3x+1)$ equals $k$ with probability $2^{-k}$, so
$\P(M=m)=2^{-(m+1)}$ for $m\ge 0$. If the $M_i$ along an orbit behave like independent random
variables with this law (\cref{lem:bernoulli} makes this precise in a 2-adic sense), then
$\log_2 x_i$ performs a random walk, and the sign of its drift, the expected number of bits gained
per odd step, decides whether orbits should come down or go up. For the Collatz map $r_i=0$ and
the drift is $\log_2 3-2\approx-0.415$: each odd step multiplies by $3$ and removes two factors $2$
on average. For $T_S$ the number $r_i$ is a function of $M_i$ determined by $S$, and the drift is
$(1+\rho(S))\log_2 3-2$, where $\rho(S)$ is the expected value of $r_i$.

\subsection*{Results}
\begin{enumerate}[leftmargin=*,label=(\arabic*)]
\item \emph{Odd-to-odd reduction and a random model} (\cref{sec:maps,sec:model}).
Between consecutive odd numbers, $T_S$ acts as $x\mapsto U_S(x)=3^{r_S(M)}\Syr(x)$, where
$\Syr(x)=(3x+1)/2^{v_2(3x+1)}$ is the Syracuse map, $M=v_2(3x+1)-1$, and $r_S(M)$ is the number of
consecutive integers $M,M-1,\dots$ lying in $S$. The valuations along a $U_S$-orbit are
2-adically independent and geometric, so $\log_2$ of an orbit behaves like a random walk with
drift $\delta(S)=(1+\rho(S))\log_2 3-2$, where $\rho(S)=\sum_{m\ge 1}2^{-(m+1)}r_S(m)$.
The model predicts that almost every orbit is eventually periodic if
$\rho(S)<\theta:=\log_3(4/3)\approx 0.26186$, and that almost every orbit diverges if
$\rho(S)>\theta$. For $S=\{a,\dots,b\}$ the first case occurs exactly when $a\ge 3$, or $a=2$ and
$b\le 3$, or $a=b=1$ (\cref{prop:intervals}). The case $\rho(S)=\theta$ never occurs, although
$\rho(S)$ comes within $2\cdot 10^{-16}$ of $\theta$ (\cref{prop:nocritical}).
\item \emph{Rigorous density results} (\cref{sec:rigorous}). We prove the analogues of the theorems
of Terras, Everett and Korec on both sides of the threshold (\cref{thm:density}), and show that for
$S\neq\N$ every odd residue class modulo $2^L$ contains a positive rational number with a periodic
orbit (\cref{prop:rational-cycles}); hence no argument that uses only residues can prove that the
orbits in a residue class diverge.
\item \emph{Numerical evidence} (\cref{sec:numerics}). For seven contracting maps $T_S$, including
the three intervals $\{1\}$, $\{2,3\}$ and $\{3,4,5,\dots\}$ with $\rho=1/4$, we checked every
$n<2^{34}$ (every $n<16{,}681{,}795{,}584$ for $S=\{1\}$): every orbit reaches one of the one to five
cycles listed in \cref{tab:cycles}. Drifts and stopping-time distributions agree with the model to
three to five digits, and delays, excursion exponents and basin proportions to within a few
percent. Three expanding maps behave as predicted as well.
\item \emph{Going down at odd numbers} (\cref{sec:minus}). The drift of $T^-_S$ is
$\delta(S)-\log_2 3=\rho(S)\log_2 3-2\le\log_2 3-2$, so every $T^-_S$ is contracting, with equality
only for the Collatz map $T^-_\N$. If $S$ contains no set $[c,\infty)$ and $(3/2)^L<2^a$ for every
maximal run $\{a,\dots,a+L-1\}$ of consecutive elements of $S$, then no odd step of $T^-_S$ increases
a number, and every orbit reaches the fixed point $1$ or the cycle $(2,3)$ (\cref{prop:simple}); this
covers, for instance, every $S$ without two consecutive elements, and the primes. For every other
$S$, a negative rational number, viewed as a 2-adic integer, repeats an increasing odd step forever,
and residue arguments cannot prove that all orbits descend. The first members beyond the criterion, such as $S=\{1,2\}$, have a much
stronger downward pull than the Collatz map; for seven of them we checked all $n<2^{34}$. There is also
a partial transfer to the Collatz map (\cref{rem:transfer-minus}).
\item \emph{Branching maps} (\cref{sec:branching}). For a family $\cS$ of sets we let a branch
apply any $T_S$, $S\in\cS$, at each step, and call $\cS$ \emph{bounded} if every $n$ reaches only
finitely many numbers. Every family is equivalent to a pair $S^-\subseteq S^+$
(\cref{lem:pairs}); reachable sets grow when $S^-$ shrinks or $S^+$ grows (\cref{prop:transfer}), so boundedness of a
family implies eventual periodicity for each member, and boundedness of $\{\emptyset,\{3\}\}$, for
example, would imply that no Collatz orbit diverges. A branching-random-walk heuristic
(\cref{sec:brw}) predicts that $\{\emptyset,\{3\}\}$ is bounded while $\{\emptyset,\{1\}\}$ and
$\{\emptyset,\{2\}\}$ are not. For $\{\emptyset,\{1\}\}$ we prove that a rule depending on
$x\bmod 32$ has positive drift $\frac{10}{17}\log_2\frac98$ almost surely
(\cref{thm:drift}), so that for almost all $n$ the reachable set contains numbers larger than
$n^{(2/3)\log_2 3-\varepsilon}$ for every $\varepsilon>0$, where $(2/3)\log_2 3=1.0566\ldots$
(\cref{thm:large}); on the other hand, every odd
residue class modulo $2^{10}$ contains a positive rational number whose reachable set is finite,
which rules out a proof of unboundedness based only on residues modulo $2^{10}$
(\cref{prop:rational-traps}).
\item \emph{Cycles} (\cref{app:cycles}). A cycle occurs exactly for the sets $S$ in an interval
determined by the exponents at which it goes up and down, so it can be named by the map of this
interval closest to the Collatz map and by its smallest element. A search over all $S$ at once finds
$53$ cycles of maps $T_S$ with smallest element below $10^5$ (and largest below $2^{120}$), and shows
that no map $T^-_S$ has cycles other than the shifted Collatz cycles $(1)$ and $(2,3)$ with smallest
element below $10^6$.
\item \emph{Unified maps} (\cref{sec:unified}). The maps $T_{\So,\Se}$ go up at an odd $x$ if
$v_2(x+1)\in\So$ and at an even $x$ if $v_2(x)\in\Se$, so that $T_S=T_{\N,S}$ and
$T^-_S=T_{\emptyset,S}$. Their orbits split into blocks, each made of steps up at odd numbers, at
most one step down at an odd number, and the steps at even numbers that follow, and the random model
extends to blocks. If $1\notin\So$, that is, if the odd $x\equiv 1\pmod 4$ go down, no block
multiplies $x+1$ by more than the corresponding block of the Collatz map $T^-_\N$, and infinitely
many such maps, among them $T_{[2,\infty),S}$ for every $S$, have exactly the block map of $T^-_\N$
(\cref{prop:down}); a census over all these maps finds only the cycles $(1)$, $(2,3)$ and $(3,5)$
with smallest element below $10^9$. If $1\in\So$, the census finds thousands of cycles, most of them
in maps close to the threshold, but single maps have about as many cycles as maps $T_S$ with the
same drift. So the maps between the two families do not seem to behave in a new way.
\end{enumerate}

\subsection*{Organization}
\Cref{sec:maps} reduces $T_S$ to a map between consecutive odd numbers and extends it to rational
numbers with odd denominator and to 2-adic integers. \Cref{sec:model} develops the random model and
the threshold $\theta$, \cref{sec:rigorous} proves the density results and the obstruction coming
from rational cycles, and \cref{sec:numerics} compares the model with computations. \Cref{sec:minus}
treats the maps $T^-_S$, \cref{sec:branching} the branching maps, \cref{sec:unified} the unified
maps, which contain both families, and \cref{sec:open} collects conjectures and open problems.
\Cref{app:cycles} explains how cycles are named and lists the cycles found by a search over all maps
at once. \Cref{sec:numerics,sec:minus,sec:branching} can be read largely independently of each
other; they use \cref{sec:maps,sec:model} and the method of proof of \cref{thm:density}.
\Cref{sec:unified} also uses \cref{sec:minus,app:cycles}.

\subsection*{Related work}
For the Collatz map itself, Terras \cite{Ter76} and Everett \cite{Eve77} proved that almost every
$n$ (in the sense of natural density) has an orbit that falls below $n$; Korec \cite{Kor94} showed
that almost every orbit falls below $n^{c}$ for every $c>\log 3/\log 4$, and Tao \cite{Tao22} that
almost all orbits (in the sense of logarithmic density) attain almost bounded values. The conjecture has been verified for all
$n<2^{71}$ \cite{Bar25}. The random-walk heuristic described above appears in Crandall's paper \cite{Cra78}; refined
stochastic models are due to Lagarias and Weiss \cite{LW92} and are surveyed in \cite{KL10}.

Generalized Collatz maps that act affinely on the residue classes modulo a fixed integer were
introduced by Matthews and Watts \cite{MW84} and studied with Markov-chain methods by Leigh
\cite{Lei86} and Matthews \cite{Mat10}; for $(d,g,h)$-maps Kontorovich and Sinai \cite{KS02} proved
that trajectories behave like random walks in a precise sense. Conway \cite{Con72} and Kurtz and
Simon \cite{KS07} showed that natural questions about such maps are undecidable. The maps $T_S$
and $T^-_S$ are of a different kind: the rule depends on the 2-adic valuation $m$ of $x$ rather
than on $x$ modulo a fixed number, and the step $x\mapsto k$ removes all factors of $2$ at once. Stochastic
models of the $3x+1$ problem are surveyed in \cite{KL10}; the ones we use for single maps are
the basic random-walk models of \cite{Cra78,LW92}. Rational numbers with odd denominator, which appear in
\cref{prop:rational-cycles,prop:rational-traps}, were used for the $3x+1$ map by Lagarias
\cite{Lag90}, who studied its rational cycles. A two-player game in which the moves are
$3n+1$ or $3n-1$ followed by halving was analysed by Alth\"ofer, Hartisch and Zipproth
\cite{AHZ23}, and Alth\"ofer and Zipproth \cite{AZ25} proved divergence of all orbits for a
Collatz-type recursion with an irrational multiplier. We have not found the maps $T_S$ and
$T^-_S$ or the branching maps of \cref{sec:branching} in the literature surveyed in
\cite{LagBib1,LagBib2,Lag10}, and the orbits above are not in the OEIS.

\subsection*{Notation and conventions}
For a nonzero integer $x$, $v_2(x)$ is the exponent of the largest power of $2$ that divides $x$,
and $v_2(p/q)=v_2(p)-v_2(q)$ for rational numbers. $\Zt$ is the ring of 2-adic integers. Readers
not familiar with it may think of an element of $\Zt$ as an infinite binary expansion
$\sum_{i\ge 0}a_i2^i$, $a_i\in\{0,1\}$, with addition and multiplication carried out with carries
as for integers. The integers are contained in $\Zt$ (for instance $-1=1+2+4+\cdots$), and so are
the rational numbers with odd denominator (for instance $-1/3=1+4+16+\cdots$, since
$3(1+4+16+\cdots)=-1$). Two elements of $\Zt$ are congruent modulo $2^k$ if their first $k$
binary digits agree, and a residue class modulo $2^k$ is the set of elements with prescribed first
$k$ digits. $\Zu=1+2\Zt$ is the set of odd elements, those with $a_0=1$; they are the units of
$\Zt$. The Haar probability measure $\mu$ on $\Zu$ is the law of a random odd 2-adic integer whose
digits $a_1,a_2,\dots$ are independent fair coin flips, so that every residue class modulo
$2^{k+1}$ inside $\Zu$ has measure $2^{-k}$. Statements that hold for $\mu$-almost every $x$ are the
2-adic counterpart of statements about a proportion tending to $1$ of the odd integers in long
intervals; the passage from one to the other is made in the proof of \cref{thm:density}.

A set $A$ of odd positive integers has (natural) density $d$ if
$\#\{n\in A: n\le N\}\,/\,\#\{n\le N: n\text{ odd}\}\to d$ as $N\to\infty$; densities are taken among
odd integers unless stated otherwise, and \emph{almost all} $n$ means all $n$ outside a set of
density zero. For $a\le b$ we write $[a,b]=\{a,a+1,\dots,b\}$ and $[a,\infty)=\{a,a+1,\dots\}$;
$\even$ and $\odd$ are the even and the odd positive integers.

\section{The maps \texorpdfstring{$T_S$}{T\_S}}\label{sec:maps}

\Cref{sec:maps,sec:model,sec:rigorous,sec:numerics} are about the maps $T_S$; the maps $T^-_S$ are
treated in \cref{sec:minus}, where most of the arguments carry over with small changes.

\subsection{The odd-to-odd map}
For odd $x$ (an odd integer, or more generally an element of $\Zu$ other than $-1/3$) put
\[
M(x)=v_2(3x+1)-1\ \ge 0,\qquad \Syr(x)=\frac{3x+1}{2^{M(x)+1}},
\]
so that $\Syr$ is the Syracuse map and $(3x+1)/2=2^{M(x)}\Syr(x)$. For $S\subseteq\N$ and
$M\ge 0$ let
\begin{equation}\label{eq:r}
r_S(M)=\max\{j\in\{0,\dots,M\}: M,M-1,\dots,M-j+1\in S\},
\end{equation}
the length of the run of consecutive elements of $S$ that ends at $M$; equivalently $r_S(0)=0$,
and $r_S(M)=r_S(M-1)+1$ if $M\in S$, $r_S(M)=0$ if $M\notin S$. Define
\begin{equation}\label{eq:US}
U_S(x)=3^{\,r_S(M(x))}\,\Syr(x)=\frac{3^{\,r_S(M)}\,(3x+1)}{2^{M+1}}.
\end{equation}

For example, $r_{\{2,3\}}(m)=0,1,2,0,0,\dots$ for $m=1,2,3,4,5,\dots$; $r_{\even}(m)$ is $1$ for
even and $0$ for odd $m$; and $r_\N(m)=m$. The next lemma shows that $U_S$ maps each odd number of a
$T_S$-orbit to the next one. After the step $x\mapsto(3x+1)/2$, the map $T_S$ trades a factor $2$
for a factor $3$ as long as the current exponent of $2$ lies in $S$, and at the first exponent
outside $S$ it removes all remaining factors $2$; so only the run of consecutive elements of $S$
that ends at $M$ matters.

\begin{lemma}\label{lem:US}
Let $x$ be odd, $M=M(x)$ and $r=r_S(M)$. The $T_S$-orbit of $x$ is
\[
x,\ \tfrac{3x+1}{2},\ \tfrac32\cdot\tfrac{3x+1}{2},\ \dots,\ \left(\tfrac32\right)^{r}\tfrac{3x+1}{2},\ \text{and then } U_S(x)
\ \text{ if } r<M,
\]
while for $r=M$ the last displayed term is already $U_S(x)$. Thus $T_S$ reaches the odd number $U_S(x)$ after exactly
\[
\ell_S(M)=1+r+[\,r<M\,]
\]
steps, and all intermediate terms are even and larger than $x$. Similarly, an even number
$2^m k$ ($k$ odd) reaches the odd number $3^{r_S(m)}k$ after $r_S(m)+[\,r_S(m)<m\,]$ steps.
\end{lemma}

\begin{proof}
$T_S(x)=(3x+1)/2=2^M\Syr(x)$. A number $3^i2^{e}\Syr(x)$ with $e\ge 1$ is multiplied by $3/2$ if
$e\in S$ and replaced by its odd part $3^i\Syr(x)$ otherwise. Starting from $e=M$, the exponent
decreases by one at each multiplication, and by \eqref{eq:r} exactly $r$ multiplications occur
before the exponent reaches $0$ (if $r=M$) or an exponent $M-r\notin S$ (if $r<M$), which is
followed by one division step. The same argument applies to $2^m k$.
\end{proof}

For instance, let $S=\{1\}$. For $x=17$ we have $3x+1=52=2^2\cdot 13$, so $M=1$, $\Syr(17)=13$ and
$r=r_S(1)=1=M$: the orbit is $17,26,39$, and $U_S(17)=39$ is reached after $\ell_S(1)=2$ steps. For
$x=453$ we have $3x+1=2^4\cdot 85$, so $M=3$ and $r_S(3)=0$: the orbit is $453,680,85$, and
$U_S(453)=\Syr(453)=85$ is again reached after $\ell_S(3)=2$ steps. Both pieces occur in the first
example orbit of \cref{sec:intro}.

Consequently the $T_S$-orbit of any $n$ is eventually periodic if and only if the $U_S$-orbit of
the first odd number in it is, and cycles of $T_S$ correspond to cycles of $U_S$. For $S=\emptyset$
we have $U_\emptyset=\Syr$, and for $S=\N$, $U_\N(x)=3^{M}\Syr(x)$.

\begin{lemma}[Monotonicity]\label{lem:monotone}
If $S\subseteq S'$ then $r_S\le r_{S'}$, $U_S\le U_{S'}$ and $T_S\le T_{S'}$ pointwise.
\end{lemma}

\begin{proof}
Immediate from \eqref{eq:r}, \eqref{eq:US} and \eqref{eq:TS}.
\end{proof}

The maps $T_S$ are not monotone in $x$ (for $3\notin S$, for instance, $7<8$ but
$T_S(7)=11>1=T_S(8)$). So the pointwise inequality $T_S\le T_{S'}$ does not imply
$T_S^{\,t}\le T_{S'}^{\,t}$, and \cref{lem:monotone} does not compare orbits; the branching maps of
\cref{sec:branching} are designed to make such comparisons.

\subsection{Rationals and 2-adic integers}\label{sec:rationals}
Extending the maps beyond the positive integers serves two purposes. The 2-adic integers carry a
natural probability measure, under which the valuations along an orbit become independent random
variables (\cref{lem:bernoulli}). The rational numbers with odd denominator, on the other hand, are
2-adic integers that are not integers: they share all residue information with integers, but
their orbits can behave differently, and they show which statements cannot be proved from residues
alone.

The formulas \eqref{eq:TS} and \eqref{eq:US} make sense for every positive rational number with odd
denominator: write $x=2^m k$ where $k$ has odd numerator and odd denominator. They also define
$U_S$ on $\Zu\setminus\{-1/3\}$. Every orbit of positive rationals with odd denominator stays in this
set, and for odd $x$ the identity
\begin{equation}\label{eq:ratio}
\frac{U_S(x)}{x}=\frac{3^{1+r_S(M)}}{2^{M+1}}\left(1+\frac{1}{3x}\right)
\end{equation}
holds for every positive $x$ with odd denominator. We use these extensions in
\cref{prop:rational-cycles,prop:rational-traps}.

The simplest rational cycles are the following. If $m\notin S$, then $y=1/(2^{m+1}-3)$ satisfies
$3y+1=2^{m+1}y$, so $T_S(y)=2^m y$ and $T_S(2^m y)=y$: the pair $(y,2^my)$ is a cycle of $T_S$. For
$m=1$ it is the cycle $(1,2)$ of the Collatz map, for $m=2$ the cycle $(1/5,4/5)$, and for $m=3$
the cycle $(1/13,8/13)$.

\section{The random model}\label{sec:model}

The following lemma is the source of all randomness in this paper.

\begin{lemma}[Independence of valuations]\label{lem:bernoulli}
Let $x\in\Zu$ be distributed according to $\mu$, and let $x_0=x$, $x_{i+1}=U_S(x_i)$ and $M_i=M(x_i)$
(these are defined for all $x$ outside a countable set). Then $M_0,M_1,\dots$ are independent, and $\P(M_i=m)=2^{-(m+1)}$ for $m\ge 0$. Moreover, for every
$m_0,\dots,m_{j-1}\ge 0$ the set of $x$ with $M_i=m_i$ for $i<j$ is a single residue class modulo
$2^{W+1}$, where $W=\sum_{i<j}(m_i+1)$, and $U_S^{\,j}$ maps it affinely onto $\Zu$.
\end{lemma}

\begin{proof}
For $m\ge 0$ let $D_m=\{x\in\Zu: M(x)=m\}$. Since $3$ is a unit,
$D_m=\{x: 3x+1\equiv 2^{m+1}\pmod{2^{m+2}}\}$ is a single residue class $a_m+2^{m+2}\Zt$ modulo
$2^{m+2}$, so $\mu(D_m)=2^{-(m+1)}$. With $r=r_S(m)$ we have
\[
U_S\bigl(a_m+2^{m+2}z\bigr)=3^{r}\,\frac{3a_m+1}{2^{m+1}}+2\cdot 3^{r+1}z\qquad(z\in\Zt),
\]
and $(3a_m+1)/2^{m+1}$ is odd. Hence $U_S$ maps $D_m$ bijectively and affinely onto $\Zu$, and it
maps the residue class modulo $2^{m+2+k}$ of $a_m+2^{m+2}z$ onto the residue class modulo $2^{k+1}$
of its image; so $\mu(U_S(E))=2^{m+1}\mu(E)$ for measurable $E\subseteq D_m$.

We prove the remaining assertions by induction on $j$, the case $j=0$ being trivial. The set $E$ of
$x$ with $M_i(x)=m_i$ for $i<j$ is $D_{m_0}\cap U_S^{-1}(E')$, where $E'$ is the set of $y$ with
$M_i(y)=m_{i+1}$ for $i<j-1$. By the induction hypothesis, $E'$ is a residue class modulo
$2^{W'+1}$, $W'=\sum_{1\le i<j}(m_i+1)$, which $U_S^{\,j-1}$ maps affinely onto $\Zu$. Its preimage
under the affine bijection $U_S\colon D_{m_0}\to\Zu$ is a residue class modulo
$2^{W'+1+m_0+1}=2^{W+1}$, which $U_S^{\,j}$ maps affinely onto $\Zu$, and its measure is
$2^{-W}=\prod_{i<j}2^{-(m_i+1)}$. Since this probability is a product of the probabilities
$2^{-(m_i+1)}$, the $M_i$ are independent with the stated law.
\end{proof}

In terms of binary digits, \cref{lem:bernoulli} says the following. If the digits of $x$ are
independent fair coin flips, then so are the digits of $U_S(x)$, and the valuations read off
consecutive disjoint blocks of digits of $x$: $M_0$ is determined by the first $M_0+2$ digits,
$M_1$ by the next $M_1+1$ digits, and so on.

The law of the valuations does not depend on $S$, because the extra factor $3^{r}$ is a 2-adic
unit. \Cref{lem:bernoulli} is the 2-adic form of Terras's observation \cite{Ter76} that the
parity sequence of the first $j$ steps of $T$ is determined by, and uniformly distributed with,
$n\bmod 2^j$. By \eqref{eq:ratio}, for a positive odd integer $x$,
\begin{equation}\label{eq:X}
\log_2 U_S(x)-\log_2 x = X_S(M)+\log_2\!\left(1+\tfrac{1}{3x}\right),\qquad
X_S(M)=(1+r_S(M))\log_2 3-(M+1).
\end{equation}
The \emph{random model} replaces the valuations along an integer orbit by an i.i.d.\ sequence
$M_i$ with $\P(M_i=m)=2^{-(m+1)}$ and ignores the small correction term; then
$\log_2 x_t-\log_2 x_0$ becomes the random walk $W_t=\sum_{i<t}X_S(M_i)$.
Here $\log_2$ measures size in bits: $X_S(M)$ is the number of bits gained (or lost, if it is
negative) in one odd step, up to the correction term, which is smaller than $1/(3x\ln 2)$.

The random model is a heuristic. For an integer $x$ all valuations along the orbit are determined
by $x$, and \cref{lem:bernoulli} shows only that, for a random odd integer from a long interval,
the first $j$ valuations are distributed like independent geometric variables as long as
$\sum_{i<j}(M_i+1)$ is smaller than the number of binary digits of the integers in question. The
rigorous results of \cref{sec:rigorous} use the model only in this range; the predictions of
\cref{sec:predictions} use it for whole orbits.

\subsection{Drift and threshold}
Since $\E(M+1)=2$, the drift of the model is
\begin{equation}\label{eq:drift}
\delta(S)=\E X_S(M)=(1+\rho(S))\log_2 3-2,\qquad \rho(S)=\E\, r_S(M)=\sum_{m\ge 1}\frac{r_S(m)}{2^{m+1}}.
\end{equation}
We call $S$ \emph{contracting} if $\delta(S)<0$ and \emph{expanding} if $\delta(S)>0$.
In words, $\delta(S)$ is the expected number of bits gained per odd step and $\rho(S)$ the expected
number of extra factors $3$ per odd step. For instance, for $S=\{1\}$ we have $r_S(M)=1$ if $M=1$,
which has probability $1/4$, and $r_S(M)=0$ otherwise; so $\rho(\{1\})=1/4$ and
$\delta(\{1\})=\frac54\log_2 3-2\approx-0.0188$.

\begin{proposition}\label{prop:threshold}
\begin{enumerate}[label=\textup{(\alph*)}]
\item $S$ is contracting if and only if $\rho(S)<\theta$, where
$\theta=\frac{2}{\log_2 3}-1=\log_3\frac43=0.261859\ldots$
\item $0=\rho(\emptyset)\le\rho(S)\le\rho(\N)=1$, and $\rho$ is monotone: $S\subseteq S'$ implies
$\rho(S)\le\rho(S')$.
\item $\theta$ is irrational. If $S$ is eventually periodic, then $\rho(S)$ is rational, so
$\delta(S)\neq 0$.
\end{enumerate}
\end{proposition}

\begin{proof}
(a) and (b) follow from \eqref{eq:drift} and \cref{lem:monotone}; for $S=\N$, $r_\N(m)=m$ and
$\rho(\N)=\sum_m m2^{-(m+1)}=1$. (c) If $\theta=p/q$ then $3^{p+q}=4^q$, which is impossible for
$q\ge 1$. If $S$ is eventually periodic, then $r_S(m)$ is either eventually periodic or, when
$S$ contains a tail $[c,\infty)$, eventually equal to $m-c+1$; in both cases the series for
$\rho(S)$ is a finite combination of geometric series and their derivatives with rational ratios.
\end{proof}

\begin{proposition}[Intervals]\label{prop:intervals}
For $S=[a,b]$ with $n=b-a+1$ elements,
$\rho(S)=2^{-a}\bigl(2-(n+2)2^{-n}\bigr)$, and $\rho([a,\infty))=2^{1-a}$. Hence $[a,b]$ is contracting
if and only if $a\ge 3$, or $a=2$ and $b\le 3$, or $a=b=1$; and $[a,\infty)$ is contracting if
and only if $a\ge 3$. The largest value of $\rho$ on contracting intervals is $1/4$, attained exactly
by $\{1\}$, $\{2,3\}$ and $[3,\infty)$.
\end{proposition}

\begin{proof}
Here $r_S(m)=m-a+1$ for $a\le m\le b$ and $0$ otherwise, so
$\rho(S)=\sum_{j=1}^{n} j\,2^{-(a+j)}=2^{-a}\bigl(2-(n+2)2^{-n}\bigr)$. For $a\ge 3$ we get
$\rho<2^{1-a}\le 1/4<\theta$. For $a=2$ the values are $1/8,1/4,11/32,\dots$ (increasing in $n$),
and $11/32>\theta$. For $a=1$ they are $1/4,1/2,\dots$ The case $b=\infty$ is the limit
$n\to\infty$.
\end{proof}

For example, $\rho(\even)=1/6$ (contracting) while $\rho(\odd)=1/3$ (expanding).

The values of $\rho$ range from $\rho(\emptyset)=0$ to $\rho(\N)=1$ and come very close to $\theta$;
the census of cycles in \cref{app:cycles}, for instance, contains the contracting set
$\{1,6,8,9\}$ with $\rho=0.26172$. Nevertheless the critical value is never attained.

\begin{proposition}[No critical sets]\label{prop:nocritical}
Let $P=\{1,6,8,9,12,15\}\cup[19,23]\cup[25,26]\cup[31,34]\cup[36,37]\cup[39,42]\cup[44,45]\cup[47,52]$,
$A=P\cup[54,\infty)$ and $B=P\cup\{53\}$. Then $\rho(A)<\theta<\rho(B)$,
$\rho(B)-\rho(A)=5\cdot 2^{-54}$, and every $S\subseteq\N$ satisfies $\rho(S)\le\rho(A)$ or
$\rho(S)\ge\rho(B)$. In particular $\rho(S)\neq\theta$, so every $S\subseteq\N$ is either contracting or
expanding.
\end{proposition}

\begin{proof}
Let $m\ge 0$, $Q\subseteq[1,m]$ and $r=r_Q(m)$. Every $S$ with $S\cap[1,m]=Q$ satisfies
$r_S(m+j)\le r+j$ for $j\ge 1$, and therefore
\[
c(Q)\le\rho(S)\le c(Q)+\frac{r+2}{2^{m+1}},\qquad c(Q)=\sum_{k=1}^{m}\frac{r_Q(k)}{2^{k+1}},
\]
with equality on the left for $S=Q$ and on the right for $S=Q\cup[m+1,\infty)$, since
$\sum_{j\ge 1}(r+j)2^{-(m+1+j)}=(r+2)2^{-(m+1)}$. Call $Q$ \emph{admissible} if this interval
contains $\theta$. A computation with exact rational arithmetic and a rigorous enclosure of
$\theta$ (ancillary file \texttt{rho\_gap\_check.py}) lists the admissible sets $Q\subseteq[1,m]$ for $m=1,2,\dots$: there are at most $12$ of them
for each $m\le 52$, and none for $m=53$. Hence for every $S$ there is a smallest $m$ such that
$S\cap[1,m]$ is not admissible, and all sets with this prefix lie on the same side of $\theta$; in
particular $\rho(S)\ne\theta$. The largest right endpoint among the intervals of these first
non-admissible prefixes that lie below $\theta$ is $\rho(A)$, attained by $A=P\cup[54,\infty)$, and
the smallest left endpoint among those above $\theta$ is $\rho(B)$, attained by $B$.
\end{proof}

The computation is small: it
handles at most twelve intervals at each of the $53$ levels and takes a fraction of a second.
Numerically $\theta-\rho(A)=1.47\cdot 10^{-16}$ and $\rho(B)-\theta=1.30\cdot 10^{-16}$. The gap
reflects a general feature of $\rho$: if a prefix $Q\subseteq[1,m]$ ends with a run of length
$r\ge 2$, the sets that contain $m+1$ have $\rho\ge c(Q)+(r+1)2^{-(m+2)}$, while those that do not
have $\rho\le c(Q)+2\cdot 2^{-(m+2)}$. So the set of values of $\rho$ is a compact set with many
gaps, and $\theta$ lies in the one that opens after the run $[47,52]$.

\begin{remark}[Single halvings]\label{rem:single}
If a downward step at an even number only halves it ($x\mapsto x/2$ for $m\ge1$, $m\notin S$), the
orbit of $(3x+1)/2=2^M\Syr(x)$ passes through the exponents $M,M-1,\dots,1$ and is multiplied by $3$
at each exponent in $S$. So $r_S(M)$ is replaced by $|S\cap[1,M]|$, and $\rho(S)$ by
$\sum_{s\in S}\P(M\ge s)=\sum_{s\in S}2^{-s}$. The results of this section and of
\cref{sec:rigorous} hold for this variant with these changes; for instance, $S=\{1\}$ is then
expanding ($\rho=1/2$) and $S=\{2\}$ contracting ($\rho=1/4$). In this paper we keep the full
division, with which the runs of consecutive elements of $S$ determine the behaviour.
\end{remark}

\subsection{Further predictions}\label{sec:predictions}
The model makes several quantitative predictions, which we test in \cref{sec:numerics}. Let
$W_t$ be the random walk above, with increments $X_i=X_S(M_i)$.
\begin{enumerate}[leftmargin=*,label=(\roman*)]
\item \emph{Steps per odd step.} By \cref{lem:US} a $T_S$-orbit spends on average
\[
\ell(S)=\E\,\ell_S(M)=1+\rho(S)+\sum_{m\ge 1,\ r_S(m)<m}2^{-(m+1)}
\]
steps per odd step.
\item \emph{Glide.} The glide of $n$ (also called its stopping time) is the number of steps of the
orbit before it first drops below $n$. In the model the glide of an odd $n$, counted in odd steps,
is $\tau=\min\{t\ge 1:W_t<0\}$; by \cref{lem:US} the orbit can drop below $n$ only at the end
of an odd step. For contracting $S$, $\E\tau<\infty$ and $\P(\tau>t)$ decays exponentially; for
expanding $S$, $p_\infty(S)=\P(\tau=\infty)>0$.
\item \emph{Delay.} For contracting $S$, the model predicts that $n$ needs about
$\ell(S)\log_2 n/|\delta(S)|$ steps to reach its cycle. For $S=\emptyset$, if each halving is counted as
a step (as for the map $T$), this is the classical estimate $2\log n/\log(4/3)\approx 6.95\log n$.
\item \emph{Excursions.} For contracting $S$ let $\lambda(S)>0$ be the root of
$\E\,2^{\lambda X}=1$. (The function $\lambda\mapsto\E\,2^{\lambda X}$ is convex, equals $1$ at
$\lambda=0$ and has derivative $\delta(S)\ln 2<0$ there; for eventually periodic $S$ it tends to
$\infty$ at the end of the interval on which it is finite, so the positive root exists and is
unique.) A random walk with negative drift climbs $h$ units above its starting point only with a
probability that decays exponentially in $h$, and Cram\'er's estimate \cite{Fel71} identifies the
rate: $\P(\sup_t W_t>h)\sim C\,2^{-\lambda(S)h}$. So the largest value on the orbit of $n$ exceeds
$n\,2^{h}$ with probability about $2^{-\lambda h}$. Consequently the record values of the orbit
maxima should grow like $n^{1+1/\lambda(S)+o(1)}$, and, for $0<\gamma<1/\lambda(S)$, the
number of $n\le N$ whose orbit exceeds $n^{1+\gamma}$ should be $N^{1-\lambda\gamma+o(1)}$.
For $S=\emptyset$, $\E\,2^{X}=3\,\E\,2^{-(M+1)}=3\sum_{j\ge 1}4^{-j}=1$, so $\lambda(\emptyset)=1$
exactly; this recovers the maximum excursion constant $2$ of the random-walk model of Lagarias and
Weiss \cite{LW92,KL10}, consistent with the computations of Oliveira e Silva \cite{OeS99}.
\item \emph{Expanding maps.} For expanding $S\neq\N$ let $\kappa(S)>0$ be the root of
$\E\,2^{-\kappa X}=1$, which exists and is unique by the same convexity argument, since $X<0$ with
positive probability when $S\neq\N$. Then $\P(\inf_t W_t<-h)\sim C\,2^{-\kappa h}$, so the proportion of
$n\approx 2^k$ whose orbit comes down to a cycle should decay like $2^{-\kappa k}=n^{-\kappa}$.
\end{enumerate}

\begin{table}[t]
\centering\small
\caption{Model constants for the ten maps studied numerically: $\rho$, drift $\delta$ (bits per
odd step), steps per odd step $\ell$, predicted delay per bit $\ell/|\delta|$, excursion exponent
$\lambda$, model mean glide $\E\tau$ (odd steps), and for expanding maps the return exponent
$\kappa$ and the never-return probability $p_\infty$. $\E\tau$ and $p_\infty$ are Monte Carlo
estimates (at least $1.5\cdot10^8$ and $2\cdot10^7$ walks); the other entries are exact or computed
to the digits shown.}
\label{tab:theory}
\begin{tabular}{lcccccc}
\toprule
$S$ & $\rho(S)$ & $\delta(S)$ & $\ell(S)$ & $\ell/|\delta|$ & $\lambda(S)$ & $\E\tau$ \\
\midrule
$\emptyset$ & $0$ & $-0.41504$ & $3/2$ & $3.61$ & $1$ & $3.493$ \\
$\{3\}$ & $1/16$ & $-0.31598$ & $25/16$ & $4.94$ & $0.8505$ & $4.095$ \\
$\{2\}$ & $1/8$ & $-0.21692$ & $13/8$ & $7.49$ & $0.4700$ & $7.651$ \\
$\even$ & $1/6$ & $-0.15088$ & $5/3$ & $11.05$ & $0.3853$ & $9.339$ \\
$[3,\infty)$ & $1/4$ & $-0.01880$ & $7/4$ & $93.10$ & $0.0980$ & $31.82$ \\
$\{2,3\}$ & $1/4$ & $-0.01880$ & $7/4$ & $93.10$ & $0.0352$ & $93.5$ \\
$\{1\}$ & $1/4$ & $-0.01880$ & $3/2$ & $79.80$ & $0.0223$ & $93.3$ \\
\midrule
 & & & & & $\kappa(S)$ & $p_\infty(S)$ \\
\midrule
$\odd$ & $1/3$ & $+0.11328$ & $19/12$ & -- & $0.1560$ & $0.1641$ \\
$[2,4]$ & $11/32$ & $+0.12979$ & $59/32$ & -- & $0.2200$ & $0.2688$ \\
$\{1,2\}$ & $1/2$ & $+0.37744$ & $13/8$ & -- & $0.3177$ & $0.4303$ \\
\bottomrule
\end{tabular}
\end{table}

\Cref{tab:theory} lists these constants. The three maps with $\rho=1/4$ have the same drift
$-0.0188$, about $4.5\%$ of the Collatz drift, but quite different excursion exponents, because
$\lambda$ depends on the whole distribution of $r_S(M)$ and not only on its mean. The column
$\ell/|\delta|$ is the predicted number of steps per bit of the starting value, and $1+1/\lambda$ is
the predicted exponent of the record peaks (\cref{sec:predictions}(iv)); both grow rapidly as
$\rho(S)$ approaches $\theta$.

\section{Rigorous results for single maps}\label{sec:rigorous}

The idea behind the following theorem is that a
typical odd integer $n$ with $k$ binary digits determines, through its residue modulo $2^k$, about
$k/2$ valuations $M_0,M_1,\dots$ of its orbit, since each odd step uses up $\E(M+1)=2$ digits on
average; by \cref{lem:bernoulli}, for most residue classes these valuations look like a sample of
independent geometric variables. By the law of large numbers the first $k/2$ odd steps therefore
change $\log_2 n$ by about $\delta k/2$, so the orbit reaches about $n^{1+\delta/2}$. What happens
afterwards is no longer controlled by the residue of $n$.

\begin{theorem}\label{thm:density}
Let $S\subseteq\N$ and $\delta=\delta(S)$.
\begin{enumerate}[label=\textup{(\alph*)}]
\item If $\delta<0$, then for every $c>1+\delta/2$ the set of odd $n$ whose $T_S$-orbit contains
a number smaller than $n^{c}$ has density one. In particular almost every orbit drops below its
starting value.
\item If $\delta>0$, then for every $c<1+\delta/2$ the set of odd $n$ whose $T_S$-orbit contains a
number larger than $n^{c}$ has density one.
\end{enumerate}
\end{theorem}

For $S=\emptyset$, part (a) is Korec's theorem \cite{Kor94}, since $1+\delta(\emptyset)/2=\log 3/\log 4$,
and its last sentence is the theorem of Terras and Everett \cite{Ter76,Eve77}. The statements
extend to all $n$. An even $n=2^m k$ reaches the odd number $3^{r_S(m)}k$ (\cref{lem:US}); for
fixed $m$, the odd $k$ for which $3^{r_S(m)}k$ is exceptional form a set of density zero, because
multiplication by $3^{r_S(m)}$ is injective and changes sizes by a bounded factor, and the
integers divisible by $2^{m+1}$ form a set of density $2^{-m-1}$.

\begin{proof}
Fix $\varepsilon\in(0,\frac14)$; it suffices to prove the statements for odd $n$ in a dyadic
interval $(N/2,N]$ as $N\to\infty$. Let $K=\lfloor(1-\varepsilon)\log_2 N\rfloor$ and
$j=\lfloor K/(2+\varepsilon)\rfloor$. Call a vector $(m_0,\dots,m_{j-1})$ \emph{good} if
$\sum_{i<j}(m_i+1)\le K$ and $\bigl|\sum_{i<j}X_S(m_i)-j\delta\bigr|\le\varepsilon j$.
By \cref{lem:bernoulli} and the law of large numbers (note $\E(M+1)=2<2+\varepsilon$), the
$\mu$-measure $p_N$ of the set $G$ of $x\in\Zu$ whose first $j$ valuations form a good vector tends
to $1$. By \cref{lem:bernoulli}, $G$ is a union of residue classes modulo $2^{K+1}$, and each odd
residue class modulo $2^{K+1}$ contains $N/2^{K+2}+O(1)$ integers of $(N/2,N]$. Since there are
$2^K=O(N^{1-\varepsilon})$ odd classes, the proportion of odd $n\in(N/2,N]$ in $G$ is
$p_N+O(N^{-\varepsilon})\to 1$.

Let $n\in G$ and $x_i=U_S^{\,i}(n)$. By \eqref{eq:X},
$\log_2 x_j=\log_2 n+\sum_{i<j}X_S(M_i)+\sum_{i<j}\log_2\bigl(1+\frac{1}{3x_i}\bigr)$, and
$j=\frac{1-\varepsilon}{2+\varepsilon}\log_2 n+O(1)$.

(a) If $x_i\ge n^c$ for all $i<j$, the last sum is at most $j/(3n^{c}\ln 2)<1$, so
$\log_2 x_j\le \log_2 n+j(\delta+\varepsilon)+1
=\bigl(1+\tfrac{(1-\varepsilon)(\delta+\varepsilon)}{2+\varepsilon}\bigr)\log_2 n+O(1)$.
For $\varepsilon$ small the coefficient is smaller than $c$, so $x_j<n^c$ for large $N$. Either way
the orbit contains a number smaller than $n^c$.

(b) The last sum is positive, so
$\log_2 x_j\ge\bigl(1+\tfrac{(1-\varepsilon)(\delta-\varepsilon)}{2+\varepsilon}\bigr)\log_2 n-O(1)$,
which exceeds $c\log_2 n$ for small $\varepsilon$ and large $N$.
\end{proof}

Tao's much stronger theorem \cite{Tao22} for $S=\emptyset$ suggests \cref{prob:tao} below.
The next proposition explains why a proof of divergence for an expanding map cannot be based on
residue classes alone.

\begin{proposition}\label{prop:rational-cycles}
Let $S\neq\N$. For every $L\ge 1$ and every odd $a$ there is a positive rational number $y$ with odd
denominator such that $y\equiv a\pmod{2^L}$ in $\Zt$ and the $T_S$-orbit of $y$ is periodic.
\end{proposition}

\begin{proof}
Choose $m_*\in\N\setminus S$, so that $r_S(m_*)=0$. By \cref{lem:bernoulli} there are
$m_0,\dots,m_{j-1}$ such that every $x\in\Zu$ whose first $j$ valuations are $m_0,\dots,m_{j-1}$
satisfies $x\equiv a\pmod{2^L}$ (take the first valuations of $a$ itself until
$\sum(m_i+1)\ge L-1$). For $s\ge 0$ consider the periodic sequence obtained by repeating the block
$(m_0,\dots,m_{j-1},m_*,\dots,m_*)$ with $s$ copies of $m_*$. By \cref{lem:bernoulli} there is
exactly one $y\in\Zu$ with this valuation sequence, and then $U_S^{\,j+s}(y)$ has the same
sequence, so $U_S^{\,j+s}(y)=y$. On the corresponding residue class, $U_S^{\,j+s}$ is an affine map
$x\mapsto(3^{A}x+c)/2^{B}$ with $A=\sum_{i<j}(1+r_S(m_i))+s$, $B=\sum_{i<j}(m_i+1)+s(m_*+1)$
and an integer $c>0$, so $y=c/(2^{B}-3^{A})$. Since $2^{m_*+1}\ge 4>3$, we have $2^B>3^A$ for
$s$ large, and then $y>0$. All points of the orbit of $y$ are positive rationals with odd
denominator, and the orbit is periodic.
\end{proof}

For example, let $S=\{1,2\}$, $a=7$ and $L=4$. The $U_S$-orbit $7\to 11\to 17\to 39$ has the
valuations $0,0,1$, which suffice since $\sum(m_i+1)=4\ge L-1$, and $m_*=3$. With $s=0$ the fixed point
is negative ($2^4<3^4$), but with $s=1$ the block $(0,0,1,3)$ gives $A=5$, $B=8$, $2^B-3^A=13$ and
$y=187/13$, whose $T_S$-orbit is the cycle
\[
\tfrac{187}{13}\to\tfrac{287}{13}\to\tfrac{437}{13}\to\tfrac{662}{13}\to\tfrac{993}{13}\to
\tfrac{1496}{13}\to\tfrac{187}{13};
\]
indeed $y\equiv 7\pmod{32}$, since the first three valuations determine $y$ modulo $2^5$.

For $S=\emptyset$ this is the classical correspondence between periodic parity sequences and
rational cycles of the $3x+1$ map; see \cite{Lag90}. The point of
\cref{prop:rational-cycles} is that the valuations along the orbit of $y$ agree with those of the
integer $a$ for as many steps as we like, so no finite amount of residue information distinguishes
the orbit of $y$ from integer orbits.

\begin{remark}
Suppose an argument shows that every positive integer in a residue class $a\bmod 2^L$ has an
unbounded $T_S$-orbit, using only the residues of the orbit points and inequalities such as
$U_S(x)\ge 3^{1+r_S(M)}x/2^{M+1}$ that hold for all positive $x$ (as in \eqref{eq:ratio}). The same
argument would apply to the rational number $y$ of \cref{prop:rational-cycles}, whose orbit is
bounded. So for every $S\neq\N$, a proof that some orbit diverges must use integrality, or some
other information that rational numbers with odd denominators do not share. For $S=\N$ divergence
is trivial.
\end{remark}

\section{Numerical results}\label{sec:numerics}

This section compares the predictions of \cref{sec:model} with computations for ten maps $T_S$:
the seven contracting maps $\emptyset$, $\{3\}$, $\{2\}$, $\even$, $[3,\infty)$, $\{2,3\}$ and $\{1\}$,
listed by increasing $\rho(S)$, and three expanding maps. The contracting maps range from the
Collatz map to the three maps with $\rho=1/4$, whose drift is only $4.5\%$ of the Collatz drift, so
that their orbits are much longer and climb much higher. The tables are meant to be read against
\cref{tab:theory}, which lists the corresponding model constants.

\subsection{Method}\label{sec:method}
For a given $S$ and bound $N$ we verify by induction on $n$ that every $n<N$ has an eventually
periodic $T_S$-orbit. For each $n$ we iterate $T_S$ until one of the following happens:
(i) the orbit drops below $n$, in which case $n$ ends in the same cycle as the smaller number
reached; (ii) a cycle is detected (by Brent's algorithm \cite{Bre80}, a cycle-detection method that needs
only a bounded amount of memory); since the orbit has not dropped
below $n$, all elements of this cycle are at least $n$; (iii) the orbit exceeds $2^{4096}$; (iv) $10^8$ steps are made. Outcomes (iii) and (iv) never
occurred. Since an orbit inside a cycle never drops below the smallest element of the cycle,
every cycle whose smallest element is below $N$ is found. For $n<2^{26}$ we also store the cycle
reached by $n$ and its \emph{delay}, the number of steps until the orbit first reaches the smallest
element of its cycle, computed recursively as $\mathrm{delay}(n)=g(n)+\mathrm{delay}(x)$, where $g(n)$
is the glide of $n$ and $x<n$ the first smaller value.

For every odd $n$ we record the glide (in steps and in odd steps), the \emph{peak} (the largest
value before the orbit drops below $n$), and $\log_2(x/n)$ at the drop. For example, under
$T_\emptyset$ the orbit of $7$ is $7,11,17,26,13,20,5,8,1$: the glide of $7$ is $6$ steps, or $4$ odd
steps, its peak is $26$, the recorded value is $\log_2(5/7)$, and the delay of $7$ is $8$. The measured drift is the
sum of $\log_2(x/n)$ divided by the total number of odd steps, over a block of starting values; by
Wald's identity $\E W_\tau=\delta(S)\,\E\tau$ for the walk stopped at the glide $\tau$, this ratio
estimates $\delta(S)$ without bias from the stopping rule. We also sample $400$ random odd
starting values with $B\in\{64,256,1024,4096\}$ bits for each map, follow them to their cycle (or,
for expanding maps, until they exceed $2^{2B+256}$) and record the delay per bit and the drift up to
the first time the orbit has at most $B/2$ bits.

The program is a Python driver with a C kernel (64-bit arithmetic with an overflow check, then
multi-limb integers), with checkpoints so that interrupted runs can be resumed, and a pure-Python
twin used for self-tests. The runs were made on free cloud notebooks with $4$ CPU cores; the speed
ranged from $3.6\cdot10^{7}$ starting values per second for $S=\emptyset$ to $1.2\cdot 10^{6}$ for
$S=\{1\}$, about $10$ hours in total. As a check, for $S=\emptyset$ the program reproduces the
known path records of the $3x+1$ problem \cite{OeS99,BarPR}; for example the largest peak for
$n<2^{34}$ is $10{,}361{,}199{,}457{,}202{,}525{,}864$, reached from $n=12{,}327{,}829{,}503$. (This is
half of the path record $20{,}722{,}398{,}914{,}405{,}051{,}728$ of the map $x\mapsto 3x+1$,
$x\mapsto x/2$, because $T_\emptyset$ combines the step $x\mapsto 3x+1$ with the following
halving.)

\subsection{Contracting maps}
We ran seven contracting maps: $\emptyset$, $\{3\}$, $\{2\}$, $\even$, and the three maps with
$\rho=1/4$. Every $n<2^{34}=17{,}179{,}869{,}184$ was checked, except for $S=\{1\}$, where the run
stopped at $n=16{,}681{,}795{,}584$; the statistics for $S=\{1\}$ below refer to the complete block
$[2^{32},2^{33})$, and those for the other maps to $[2^{33},2^{34})$.

\begin{table}[!htbp]
\centering\small
\caption{Cycles of the contracting maps, named as in \cref{app:cycles} ($C^+_U(n)$ is the cycle of
$T_U$ whose smallest element is $n$): length in steps, number of odd steps, and the share (in
percent) of starting values that end in the cycle, among all $n\in[2^{25},2^{26})$ and among $400$
random odd $4096$-bit numbers. No other cycle has smallest element below $2^{34}$ (below
$1.668\cdot 10^{10}$ for $S=\{1\}$).}
\label{tab:cycles}
\begin{tabular}{llrrrr}
\toprule
 & & & & \multicolumn{2}{c}{share of starting values}\\
\cmidrule(lr){5-6}
$S$ & cycle & length & odd steps & $n\in[2^{25},2^{26})$ & $4096$-bit\\
\midrule
$\emptyset$ & $C^+_\emptyset(1)$ & $2$ & $1$ & $100$ & $100$\\
\midrule
$\{3\}$ & $C^+_\emptyset(1)$ & $2$ & $1$ & $30.5$ & $27.8$\\
 & $C^+_{\{3\}}(3)$ & $4$ & $2$ & $26.9$ & $26.0$\\
 & $C^+_{\{3\}}(103)$ & $39$ & $27$ & $42.5$ & $46.3$\\
\midrule
$\{2\}$ & $C^+_\emptyset(1)$ & $2$ & $1$ & $52.4$ & $55.5$\\
 & $C^+_{\{2\}}(19)$ & $8$ & $4$ & $22.8$ & $23.3$\\
 & $C^+_{\{2\}}(47)$ & $24$ & $16$ & $24.7$ & $21.3$\\
\midrule
$\even$ & $C^+_\emptyset(1)$ & $2$ & $1$ & $1.5$ & $0.5$\\
 & $C^+_{\{4\}}(15)$ & $6$ & $4$ & $64.5$ & $67.8$\\
 & $C^+_{\{2\}}(19)$ & $8$ & $4$ & $34.0$ & $31.8$\\
\midrule
$[3,\infty)$ & $C^+_\emptyset(1)$ & $2$ & $1$ & $0$ & $0$\\
 & $C^+_{\{3\}}(3)$ & $4$ & $2$ & $1.6$ & $1.8$\\
 & $C^+_{[3,7]}(347)$ & $55$ & $26$ & $63.5$ & $65.3$\\
 & $C^+_{[3,7]}(255)$ & $126$ & $71$ & $25.5$ & $25.0$\\
 & $C^+_{[3,9]}(271)$ & $558$ & $304$ & $9.4$ & $8.0$\\
\midrule
$\{2,3\}$ & $C^+_\emptyset(1)$ & $2$ & $1$ & $47.4$ & $51.3$\\
 & $C^+_{\{2,3\}}(5)$ & $18$ & $10$ & $17.8$ & $16.8$\\
 & $C^+_{\{2\}}(19)$ & $8$ & $4$ & $6.7$ & $6.8$\\
 & $C^+_{\{2\}}(47)$ & $24$ & $16$ & $28.1$ & $25.3$\\
\midrule
$\{1\}$ & $C^+_{\{1\}}(1)$ & $5$ & $3$ & $100$ & $100$\\
\bottomrule
\end{tabular}
\end{table}

\Cref{tab:cycles} lists the cycles, named as in \cref{app:cycles}: $C^+_U(n)$ is the cycle with
smallest element $n$ of the map $T_U$, the map closest to the Collatz map in which the cycle occurs.
The $20$ cycles listed for the seven maps are $11$ different cycles; for example, the cycle
$(15,23,35,53,80,120)$ of $T_{\even}$ is $C^+_{\{4\}}(15)$, which occurs for every $S$ with $4\in S$ and
$3\notin S$. Every cycle has smallest element at most $347$, and no new cycle appears beyond that, up
to $2^{34}$. For $S=[3,\infty)$ the cycle $C^+_\emptyset(1)=(1,2)$ is reached only from $1$, $2$ and
$4$. The basin shares in $[2^{25},2^{26})$ and for random $4096$-bit numbers agree
within the sampling error (about $2.5$ percentage points), suggesting that each basin has a natural
density (\cref{conj:basins}).

\begin{table}[!htbp]
\centering\small
\caption{Model versus data for the contracting maps. Drift in bits per odd step, measured over all
odd $n$ in the top block (see text) and over $400$ random $4096$-bit starts (standard error in
units of the last digit); delay per bit for the $4096$-bit starts; mean glide in odd steps;
$\lambda$ fitted to the tail of the peak excess $\log_2(\text{peak}/n)$ for odd
$n\in[2^{32},2^{33})$. Model values: $\P(\text{glide}>64)$ computed exactly by dynamic programming
over the walk, the other glide statistics by Monte Carlo (\cref{tab:theory}).}
\label{tab:drift}
\begin{tabular}{lcccccc}
\toprule
 & \multicolumn{2}{c}{drift: data} & \multicolumn{2}{c}{delay per bit} & \multicolumn{2}{c}{$\lambda$}\\
\cmidrule(lr){2-3}\cmidrule(lr){4-5}\cmidrule(lr){6-7}
$S$ & all odd $n$ & $4096$-bit & data & model & fit & model\\
\midrule
$\emptyset$ & $-0.41503$ & $-0.4151(9)$ & $3.62$ & $3.61$ & $1.001$ & $1$\\
$\{3\}$ & $-0.31598$ & $-0.3152(8)$ & $4.95$ & $4.94$ & $0.852$ & $0.851$\\
$\{2\}$ & $-0.21692$ & $-0.2171(7)$ & $7.49$ & $7.49$ & $0.471$ & $0.470$\\
$\even$ & $-0.15088$ & $-0.1499(5)$ & $11.09$ & $11.05$ & $0.387$ & $0.385$\\
$[3,\infty)$ & $-0.01879$ & $-0.0188(1)$ & $92.5$ & $93.1$ & $0.0997$ & $0.0980$\\
$\{2,3\}$ & $-0.01880$ & $-0.0187(2)$ & $92.9$ & $93.1$ & $0.0357$ & $0.0352$\\
$\{1\}$ & $-0.01879$ & $-0.0188(2)$ & $80.4$ & $79.8$ & $0.0227$ & $0.0223$\\
\bottomrule
\end{tabular}

\medskip
\begin{tabular}{lcccccc}
\toprule
 & \multicolumn{2}{c}{mean glide (odd steps)} & \multicolumn{2}{c}{$\P(\text{glide}>64\text{ odd steps})$} & \multicolumn{2}{c}{$\P(\text{glide}>512\text{ odd steps})$}\\
\cmidrule(lr){2-3}\cmidrule(lr){4-5}\cmidrule(lr){6-7}
$S$ & data & model & data & model & data & model\\
\midrule
$\emptyset$ & $3.4927$ & $3.4926$ & $0.0004482$ & $0.0004481$ & $0$ & $0$\\
$\{3\}$ & $4.0946$ & $4.0946$ & $0.0018994$ & $0.0018997$ & $0$ & $0$\\
$\{2\}$ & $7.6511$ & $7.6512$ & $0.0150417$ & $0.0150413$ & $1.56\cdot10^{-6}$ & $1.6\cdot10^{-6}$\\
$\even$ & $9.3377$ & $9.3388$ & $0.0258428$ & $0.0258422$ & $4.73\cdot10^{-5}$ & $4.71\cdot10^{-5}$\\
$[3,\infty)$ & $31.810$ & $31.817$ & $0.0615554$ & $0.0615537$ & $0.012577$ & $0.012579$\\
$\{2,3\}$ & $93.54$ & $93.53$ & $0.1215696$ & $0.1215735$ & $0.031813$ & $0.031816$\\
$\{1\}$ & $93.18$ & $93.28$ & $0.0990913$ & $0.0990879$ & $0.027674$ & $0.027680$\\
\bottomrule
\end{tabular}
\end{table}

\begin{figure}[!htbp]
\centering
\includegraphics[width=\textwidth]{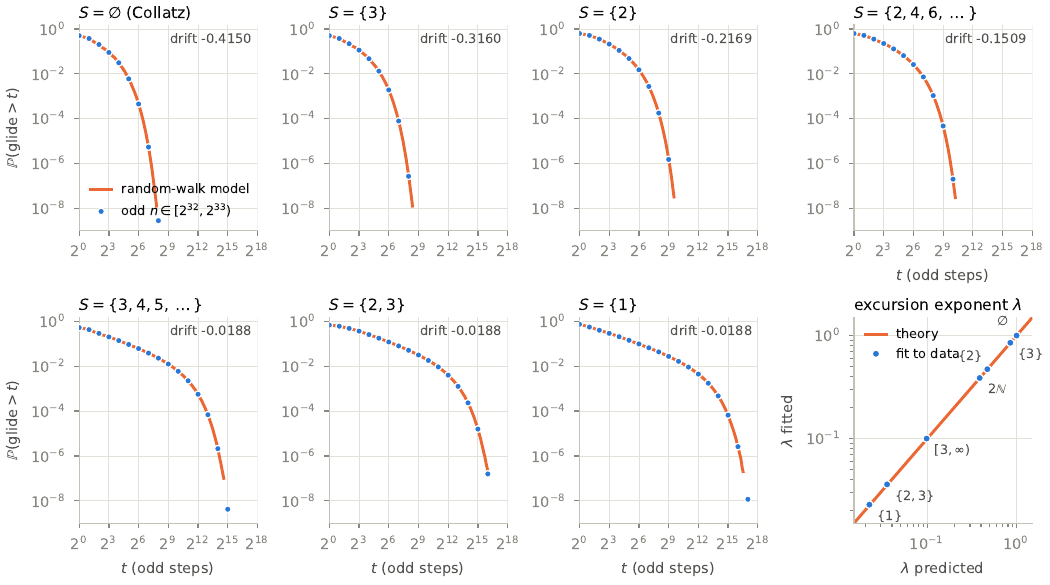}
\caption{Glide length distributions. Dots: proportion of odd $n\in[2^{32},2^{33})$ whose glide
has more than $t$ odd steps. Lines: the random-walk model of \cref{sec:model} (Monte Carlo with
$1.5\cdot10^8$ to $2\cdot10^9$ walks, drawn while at least $25$ walks remain). Last panel: the
exponent $\lambda$ fitted to the tail of the peak excess, against the root of $\E\,2^{\lambda X}=1$.}
\label{fig:glides}
\end{figure}

\Cref{tab:drift} and \cref{fig:glides} compare the data with the model. The measured drifts agree
with $\delta(S)$ to four or five digits. For every map, the proportion of glides longer than $64$
odd steps differs from the exact model value by less than one binomial standard error for
the same number of independent samples ($|z|<0.8$ in all seven cases); the mean glides and the far tails
agree to about $0.1\%$, and the fitted excursion exponents are within $2\%$ of $\lambda(S)$. The
starting values have only $33$ or $34$ bits, whereas a glide of $64$ odd steps involves about $128$
binary digits of valuations and a glide of $512$ odd steps about $1000$; so the independence of
\cref{lem:bernoulli} persists, in distribution and to high accuracy, far beyond the range where it
is a theorem.

\begin{table}[!htbp]
\centering\small
\caption{Records among $n<2^{34}$ (among $n<1.668\cdot10^{10}$ for $S=\{1\}$): the longest glide (in
steps) and the highest peak (in bits).}
\label{tab:records}
\begin{tabular}{lrrrr}
\toprule
$S$ & longest glide & at $n=$ & highest peak & at $n=$\\
\midrule
$\emptyset$ & $471$ & $12{,}235{,}060{,}455$ & $64$ & $12{,}327{,}829{,}503$\\
$\{3\}$ & $712$ & $5{,}899{,}290{,}535$ & $68$ & $15{,}886{,}636{,}263$\\
$\{2\}$ & $1{,}605$ & $15{,}860{,}082{,}535$ & $93$ & $7{,}302{,}099{,}679$\\
$\even$ & $2{,}641$ & $10{,}356{,}767{,}463$ & $110$ & $16{,}986{,}655{,}675$\\
$[3,\infty)$ & $71{,}511$ & $14{,}078{,}758{,}975$ & $290$ & $15{,}263{,}581{,}087$\\
$\{2,3\}$ & $179{,}345$ & $16{,}965{,}229{,}624$ & $746$ & $16{,}038{,}017{,}103$\\
$\{1\}$ & $251{,}984$ & $6{,}238{,}794{,}859$ & $1{,}252$ & $4{,}173{,}232{,}594$\\
\bottomrule
\end{tabular}
\end{table}

\Cref{tab:records} shows the records. The near-critical maps have glides of $10^5$ steps and peaks
of more than $1000$ bits starting from $33$-bit numbers. By \cref{sec:predictions}(iv) the ratio
$\log(\text{peak})/\log n$ of the record holders should tend to $1+1/\lambda(S)$; the observed
ratios ($1.88$, $1.99$, $2.81$, $3.22$, $8.57$, $22.0$, $39.2$ in the order of the table) are all
somewhat below the limits $2$, $2.18$, $3.13$, $3.60$, $11.2$, $29.4$, $45.8$, as expected at this
size. The three maps with $\rho=1/4$ have the same drift but differ by an order of magnitude in
their records, in the order predicted by $\lambda$.

\subsection{Expanding maps}
For the expanding maps $\odd$, $[2,4]$ and $\{1,2\}$ we checked all $n<2^{20}$, treating an orbit
as escaping when it exceeds $2^{256}$. The cycles found are $C^+_{\{3\}}(3)=(3,5,8,12)$ for $\odd$;
$C^+_\emptyset(1)$, $C^+_{\{2\}}(19)$ and $C^+_{\{2,3,4\}}(23)$ for $[2,4]$; and
$C^+_{\{1\}}(1)=(1,2,3,5,8)$ for $\{1,2\}$. The proportion of odd
$n\in[2^{19},2^{20})$ that never drop below $n$ is $0.1644$, $0.2686$ and $0.4299$, against the model
values $p_\infty=0.1641$, $0.2688$, $0.4303$. The proportion of all $n\in[2^{k-1},2^{k})$ whose orbit
reaches a cycle decreases from $19.9\%$, $9.3\%$, $8.1\%$ at $k=15$ to $11.6\%$, $4.3\%$, $2.6\%$ at
$k=20$; the least-squares slopes of its $\log_2$ over $k=15,\dots,20$ are $-0.156$, $-0.217$,
$-0.321$, against $-\kappa=-0.156$, $-0.220$, $-0.318$. All $400$ random $4096$-bit starting values
of each map escaped, with measured drifts $+0.1130(4)$, $+0.1302(3)$, $+0.3781(7)$ against
$\delta=+0.11328$, $+0.12979$, $+0.37744$.

\let\SavedFloatBarrier\FloatBarrier\let\FloatBarrier\relax
\section{Going down at odd numbers: the maps \texorpdfstring{$T^-_S$}{T-\_S}}\label{sec:minus}
\let\FloatBarrier\SavedFloatBarrier

For odd $x$, an odd integer or an element of $\Zu\setminus\{-1\}$, put $M^-(x)=v_2(x+1)-1$.
Exactly as in \cref{lem:US}, the $T^-_S$-orbit of an odd $x$ reaches the odd number
\begin{equation}\label{eq:Uminus}
U^-_S(x)=3^{\,r_S(M)}\,\frac{x+1}{2^{M+1}},\qquad M=M^-(x),
\end{equation}
after $\ell_S(M)$ steps, all intermediate terms being even. The set $\{x\in\Zu:M^-(x)=m\}$ is again a
residue class modulo $2^{m+2}$ that $U^-_S$ maps affinely onto $\Zu$, so \cref{lem:bernoulli} holds
verbatim: the valuations $M^-$ along a $U^-_S$-orbit are independent with $\P(M^-=m)=2^{-(m+1)}$.
Since $\log_2 U^-_S(x)-\log_2 x=r_S(M)\log_2 3-(M+1)+\log_2(1+1/x)$, the drift of the model is
\begin{equation}\label{eq:drift-minus}
\delta^-(S)=\rho(S)\log_2 3-2=\delta(S)-\log_2 3\ \le\ \log_2 3-2=-0.41504\ldots,
\end{equation}
with equality only for $S=\N$ (if $m_0\notin S$, then $r_S(m_0)=0<m_0$ and $\rho(S)<1$). So every $S$
is contracting for $T^-$, and the proof of \cref{thm:density}(a) applies with only notational changes
(the correction term is now $\log_2(1+1/x)$): for every
$S$ and every $c>1+\delta^-(S)/2$, almost every odd $n$ has a $T^-_S$-orbit that contains a number
smaller than $n^c$.

\subsection{The Collatz map at the top}
The maps $T^-_S$ look rather different from the Collatz map, but the largest of them is the
Collatz map in disguise.

\begin{proposition}\label{prop:shift}
For every $x\ge 1$, $T^-_S(x)=G_S(x-1)+1$, where $G_S$ is the map of the nonnegative integers
\[
G_S(y)=\begin{cases}
y/2, & y \text{ even},\\
(3y+1)/2, & y \text{ odd and } m\in S,\\
\lfloor y/2^{m}\rfloor, & y \text{ odd and } m\notin S,
\end{cases}
\qquad m=v_2(y+1).
\]
In particular $G_\N$ is the Collatz map $T$, and $T^-_\N(x)=T(x-1)+1$.
\end{proposition}

\begin{proof}
If $x$ is odd, then $y=x-1$ is even and $T^-_S(x)=(x+1)/2=y/2+1$. If $x=2^m k$ with $m\ge1$ and $k$
odd, then $y=x-1$ is odd with $v_2(y+1)=m$; if $m\in S$, then $T^-_S(x)=3x/2=(3y+1)/2+1$, and
otherwise $T^-_S(x)=k=(y+1)/2^m=\lfloor y/2^m\rfloor+1$.
\end{proof}

For odd $y$, $m=v_2(y+1)$ is the length of the final block of ones in the binary expansion of $y$,
and $\lfloor y/2^m\rfloor$ deletes this block. For example, $y=23=10111_2$ ends in a
block of three ones, so $m=3$; if $3\in S$ then $G_S(23)=(3\cdot 23+1)/2=35$, and otherwise
$G_S(23)=\lfloor 23/8\rfloor=2=10_2$. Thus $T^-_S$ is, up to the shift by one, the Collatz
map in which an odd number ending in a block of exactly $m$ ones is sent to $(3y+1)/2$ only when
$m\in S$, and otherwise loses that block. The fixed point $1$ and the cycle $(2,3)$ of $T^-_\N$ are
the Collatz cycles $(0)$ and $(1,2)$, and the Collatz conjecture is equivalent to the statement that
every $T^-_\N$-orbit of a number $x\ge 2$ reaches $(2,3)$.

\begin{remark}[A partial transfer]\label{rem:transfer-minus}
$G_{[1,b]}$ agrees with $T$ except at the odd numbers that end in at least $b+1$ ones, that is, at
$y\equiv -1\pmod{2^{b+1}}$. Hence a Collatz orbit that never meets such a number is, after the shift,
an orbit of $T^-_{[1,b]}$. Consequently, if for every $b$ every orbit of $T^-_{[1,b]}$ is eventually
periodic, then every divergent Collatz orbit contains numbers that end in arbitrarily long blocks of
ones.
\end{remark}

\subsection{An elementary criterion}\label{sec:criterion}

The simplest situation is that no odd step can increase a number. By \eqref{eq:Uminus}, an odd step
maps $x$ to $U^-_S(x)=\bigl(3^{r_S(M)}/2^{M+1}\bigr)(x+1)$, and the following criterion asks that all
these factors be smaller than $1$.

\needspace{9\baselineskip}
\begin{proposition}\label{prop:simple}
For $S\subseteq\N$ the following are equivalent:
\begin{enumerate}[label=\textup{(\roman*)}]
\item $3^{r_S(M)}<2^{M+1}$ for every $M\ge 0$;
\item no odd step of $T^-_S$ increases a positive integer, that is, $U^-_S(x)\le x$ for every odd
$x\ge 1$;
\item $S$ contains no set $[c,\infty)$, and $(3/2)^L<2^a$ for every maximal run $[a,a+L-1]\subseteq S$
of consecutive integers.
\end{enumerate}
If they hold, then $U^-_S(x)<x$ for every odd $x>3$, and every $T^-_S$-orbit of a positive integer
reaches the fixed point $1$ or, if $1\in S$, the cycle $(2,3)$. They hold in particular if
$q(S)=\sup_{M\ge 0}3^{r_S(M)}/2^{M+1}<1$.
\end{proposition}

\begin{proof}
For odd $x$ write $x+1=2^{M+1}k$ with $M=M^-(x)$ and $k$ odd, and let $r=r_S(M)$; then
$U^-_S(x)=3^rk$. If (i) holds, then $3^rk<2^{M+1}k=x+1$, which gives (ii). If (i) fails for some $M$,
then $3^r>2^{M+1}$, and $x=2^{M+1}-1$ satisfies $M^-(x)=M$ and $U^-_S(x)=3^r>x$; so (ii) implies
(i). For the equivalence of (i) and (iii): $r_S(M)=0$ for $M\notin S$, which gives the factor
$2^{-(M+1)}\le 1/2$; along a run $[a,a+L-1]$ the factor $3^{M-a+1}/2^{M+1}$ increases with $M$ and is
largest at the end of the run, where it equals $(3/2)^L2^{-a}$; and if $[c,\infty)\subseteq S$, then
$r_S(M)\ge M-c+1$, so the factor tends to infinity. Now assume (i). Then $U^-_S(x)=x$ means
$k(2^{M+1}-3^r)=1$, so $k=1$ and $2^{M+1}-1=3^r$. For $M\ge 2$ this is impossible, since
$2^{M+1}-1\equiv 7\pmod 8$ while $3^r\equiv 1$ or $3\pmod 8$. So either $M=r=0$ and $x=1$, or $M=r=1$,
$1\in S$ and $x=3$. Hence the odd elements of every orbit decrease strictly until they reach $1$ or
$3$; and $U^-_S(3)=3^{r_S(1)}$ is $3$ if $1\in S$, when $3\mapsto 2\mapsto 3$, and $1$ otherwise.
Finally, $q(S)<1$ clearly implies (i).
\end{proof}

By (iii), a run starting at $a$ is allowed if and only if its length is at most
$\lfloor a/\log_2(3/2)\rfloor=\lfloor 1.7095\ldots a\rfloor$. The criterion therefore holds, for
instance, for $S=\emptyset$ and $S=\{1\}$; for every $S$ without two consecutive elements, such as the
even numbers, the odd numbers and the squares; for the primes, whose only run of length two,
$\{2,3\}$, gives $(3/2)^2<2^2$; for the longest allowed runs starting at $a=2,3,4,10$, which are
$\{2,3,4\}$, $[3,7]$, $[4,9]$ and $[10,26]$; and for every set made of allowed runs with gaps between
them, such as $\{2,3,4\}\cup[6,15]\cup\{20\}$. It fails for $\{1,2\}$, for every $[1,b]$ with $b\ge 2$,
for $\{2,3,4,5\}$, and for every $S$ that contains a set $[c,\infty)$, in particular for $\N$.

When the criterion fails, some odd step increases $x$, and a 2-adic integer repeats that step
forever. Indeed, if $3^r>2^{M+1}$ for $r=r_S(M)$, then the step $x\mapsto 3^r(x+1)/2^{M+1}$, taken on the
class $M^-(x)=M$, has the fixed point $x^*=-3^r/(3^r-2^{M+1})$, a negative rational number with odd
denominator and $M^-(x^*)=M$. For $S=\{1,2\}$ and $M=2$ this is the step $x\mapsto 9(x+1)/8$, taken at
$x\equiv 7\pmod{16}$, with $x^*=-9$; for $S=\{2,3,4,5\}$ and $M=5$ it is the step $x\mapsto 81(x+1)/64$,
with $x^*=-81/17$. Integers close to $x^*$ imitate it: under $T^-_{\{1,2\}}$ the first nine odd steps
of $2^{30}-9$ all take the step at $M=2$ and lead to $9^j2^{30-3j}-9$, $j\le 9$, so that the number
grows from $30$ to $31.5$ bits. Recall that on a residue class modulo $2^K$ that determines the first
$j$ valuations, the $j$-th iterate of $U^-_S$ is affine with multiplier $3^a/2^b$; following Terras
\cite{Ter76}, the class has a \emph{coefficient stopping time} if $3^a/2^b<1$ for some such $j$, and
then all sufficiently large elements of the class eventually decrease. Along the residue classes that
contain $x^*$ every determined step has multiplier $3^r/2^{M+1}>1$, so no residue class modulo $2^K$
that contains $x^*$ has a coefficient stopping time, for any $K$, just as for the Collatz map and the
point $-1$. Under the criterion, on the other hand, already the first step has multiplier less than
$1$ on every class. So the criterion holds exactly when every odd 2-adic integer lies in a residue
class with a coefficient stopping time; when it fails, the descent of all orbits cannot be proved
from residue classes alone. These maps form the part of the family between the elementary cases and
the Collatz map.

\subsection{Numerical results}
We checked, with the method of \cref{sec:method}, that for each of the seven sets $S$ in
\cref{tab:minus} every $n<2^{34}$ has an eventually periodic $T^-_S$-orbit, and that the only cycles
with smallest element below $2^{34}$ are the fixed point $C^-_\emptyset(1)=(1)$ and, when $1\in S$, the cycle $C^-_\emptyset(2)=(2,3)$. An odd starting value drops at the
first step, so long glides and high peaks come from even starting values $2^m k$ with $m$ in a long
run of $S$; even so, no glide is longer than $279$ steps, and no orbit gets above $2^{54}$. The drift measured on random
$4096$-bit starting values agrees with $\delta^-(S)$ to two or three digits.

\begin{table}[!htbp]
\centering\small
\caption{Maps $T^-_S$ checked for all $n<2^{34}$: $\rho(S)$, the drift $\delta^-(S)$ and the drift
measured on $100$ random $4096$-bit starting values (followed until they have at most $2048$ bits),
the quantity $q(S)$ of \cref{prop:simple}, the longest glide (in steps) and the highest peak (in
bits). The only cycles with smallest element below $2^{34}$ are $C^-_\emptyset(1)=(1)$ and, for
$1\in S$, $C^-_\emptyset(2)=(2,3)$. The
last row is the Collatz map, for comparison.}
\label{tab:minus}
\begin{tabular}{lcccccc}
\toprule
$S$ & $\rho(S)$ & $\delta^-(S)$ & drift, samples & $q(S)$ & longest glide & highest peak\\
\midrule
$\{1,2\}$ & $1/2$ & $-1.20752$ & $-1.2098$ & $9/8$ & $45$ & $36.9$\\
$[1,3]$ & $11/16$ & $-0.91034$ & $-0.9098$ & $27/16$ & $100$ & $41.4$\\
$[1,6]$ & $15/16$ & $-0.51410$ & $-0.5107$ & $729/128$ & $279$ & $53.9$\\
$\N\setminus\{1\}$ & $1/2$ & $-1.20752$ & $-1.2097$ & $\infty$ & $94$ & $51.8$\\
$\N\setminus\{2\}$ & $1/2$ & $-1.20752$ & $-1.2034$ & $\infty$ & $82$ & $51.1$\\
$\N\setminus\{3\}$ & $5/8$ & $-1.00940$ & $-1.0105$ & $\infty$ & $80$ & $50.6$\\
$[3,\infty)$ & $1/4$ & $-1.60376$ & $-1.6027$ & $\infty$ & $58$ & $51.1$\\
\midrule
$\N$ & $1$ & $-0.41504$ & $-0.4147$ & $\infty$ & -- & --\\
\bottomrule
\end{tabular}
\end{table}

Beyond these seven maps, a search over all sets $S$ at once (\cref{app:cycles}) shows that for every
$S\subseteq\N$ the only cycles of $T^-_S$ with smallest element below $10^6$ are $C^-_\emptyset(1)$ and
$C^-_\emptyset(2)$. So, unlike the maps $T_S$ (\cref{tab:cycles,tab:census}), the maps $T^-_S$ seem
to have no cycles besides the shifted cycles of the Collatz map.

\section{Branching maps}\label{sec:branching}

\Cref{lem:monotone} compares single maps only pointwise and says nothing about their orbits. To
compare orbits we let the map vary along the orbit. One may think of a player who chooses at each
step which member of a family to apply (compare the Collatz game of \cite{AHZ23}); a family is
bounded if no strategy drives any starting value to infinity. Since every orbit of every member is
one possible play, boundedness is a statement about all members at once, and it passes from a
family to smaller families.

\subsection{Definitions}
\begin{definition}\label{def:branching}
A \emph{branching family} is a finite nonempty set $\cS$ of subsets of $\N$. For $x\in\N$ let
$N_\cS(x)=\{T_S(x):S\in\cS\}$. A \emph{branch} from $n$ is a sequence $n=x_0,x_1,x_2,\dots$
with $x_{t+1}\in N_\cS(x_t)$, and the \emph{reachable set} $R_\cS(n)$ is the set of all numbers
on branches from $n$. The family is \emph{bounded} if $R_\cS(n)$ is finite for every $n$, and
\emph{unbounded} otherwise.
\end{definition}

The member applied may change at every step, also in the middle of a run of multiplications by
$3/2$. Each $T_S$-orbit with $S\in\cS$ is a branch, so in a bounded family every member has only
eventually periodic orbits. The definition extends verbatim to positive rationals with odd
denominators. One could also allow a choice at odd numbers, between $(3x+1)/2$ and $(x+1)/2$, by
admitting maps $T^-_S$ as members; we do not consider such families here.

\begin{lemma}\label{lem:pairs}
Let $S^-=\bigcap_{S\in\cS}S$ and $S^+=\bigcup_{S\in\cS}S$. For odd $x$,
$N_\cS(x)=\{(3x+1)/2\}$; for $x=2^m k$ with $m\ge 1$, $N_\cS(x)$ contains $3x/2$ if and only if
$m\in S^+$, and contains $k$ if and only if $m\notin S^-$. Consequently $R_\cS=R_{\{S^-,S^+\}}$.
\end{lemma}

\begin{proof}
$T_S(x)=3x/2$ for some $S\in\cS$ iff $m\in S^+$, and $T_S(x)=k$ for some $S\in\cS$ iff $m\notin S^-$.
\end{proof}

We therefore write $\pair{S^-}{S^+}$ for the family $\{S^-,S^+\}$, $S^-\subseteq S^+$: at an exponent
in $S^-$ the multiplication is forced, at an exponent in $S^+\setminus S^-$ there is a choice, and
outside $S^+$ the division is forced. A single map is $\pair{S}{S}$, and for example
$\{\emptyset,\{2\},\{3\}\}$ is equivalent to $\pair{\emptyset}{\{2,3\}}$. As in \cref{lem:US}, the odd
numbers reachable from an odd $x$ in one odd step are $3^{j}\Syr(x)$ with $j$ in
\begin{equation}\label{eq:J}
J(M)=\bigl\{j\in\{0,\dots,M\}:\ M-i\in S^+\text{ for } 0\le i<j,\ \text{and } j=M \text{ or } M-j\notin S^-\bigr\},
\quad M=M(x).
\end{equation}
For example, in $\pair{\emptyset}{\{1\}}$ the odd number $x=17$ has $M=1$ and $J(1)=\{0,1\}$, so both
$13=\Syr(17)$ and $39=3\Syr(17)$ are reachable from $17$ in one odd step.

\begin{lemma}\label{lem:konig}
$R_\cS(n)$ is infinite if and only if there is a branch from $n$ consisting of pairwise distinct
numbers; such a branch tends to infinity.
\end{lemma}

\begin{proof}
Each number has at most two successors. If $R_\cS(n)$ is infinite, the tree obtained by attaching
each element of $R_\cS(n)\setminus\{n\}$ to the element from which a breadth-first search first
reaches it is infinite and locally finite, so by K\"onig's lemma it has an infinite path from $n$,
which is a branch with distinct terms. A sequence of distinct positive integers tends to infinity.
The converse is clear.
\end{proof}

\subsection{Transfer}
\begin{proposition}[Transfer principle]\label{prop:transfer}
If $A'\subseteq A\subseteq B\subseteq B'$, then $R_{\pair{A}{B}}(n)\subseteq R_{\pair{A'}{B'}}(n)$ for
every $n$. Hence, if $\pair{A'}{B'}$ is bounded, then so is $\pair{A}{B}$, and every orbit of
every $T_S$ with $A'\subseteq S\subseteq B'$ is eventually periodic. Conversely, if some $T_S$ with
$A'\subseteq S\subseteq B'$ has an unbounded orbit, then $\pair{A'}{B'}$ is unbounded.
\end{proposition}

\begin{proof}
By \cref{lem:pairs}, $N_{\pair{A}{B}}(x)\subseteq N_{\pair{A'}{B'}}(x)$ for every $x$.
\end{proof}

\begin{corollary}\label{cor:collatz}
If the family $\pair{\emptyset}{\{3\}}$ is bounded, then no orbit of the Collatz map diverges, and no
orbit of $T_{\{3\}}$ diverges.
\end{corollary}

Thus the branching families form a partially ordered set in which boundedness propagates
downwards and unboundedness upwards, and single maps are the minimal elements. A proof that one family
$\pair{S^-}{S^+}$ is bounded would settle the eventual periodicity of all maps $T_S$ with
$S^-\subseteq S\subseteq S^+$ simultaneously; when $S^+\setminus S^-$ is infinite, these are
uncountably many maps. A statement about a
family is in general stronger than the corresponding statements about its members; the families
of interest are those that are still plausibly bounded, or those that are unbounded although all
their members are contracting.

\subsection{A branching-random-walk heuristic}\label{sec:brw}
Recall that a branching random walk is a system of particles on the real line in which every
particle of generation $t$, independently of the others, is replaced in generation $t+1$ by a random
finite set of children, whose displacements relative to the parent have a fixed joint law. In the
random model the odd numbers reachable from $x$ form such a walk: a particle
at position $\log_2 x$ has children at $\log_2 x+X_j(M)$, $j\in J(M)$, where
$X_j(M)=(1+j)\log_2 3-(M+1)$ and $M$ is geometric as before. Let
\[
\Lambda(\lambda)=\E\sum_{j\in J(M)}2^{\lambda X_j(M)}
=\sum_{M\ge 0}2^{-(M+1)}\sum_{j\in J(M)}3^{\lambda(1+j)}2^{-\lambda(M+1)}.
\]
Thus $\Lambda(\lambda)$ is the expected value of $\sum 2^{\lambda\cdot(\text{displacement})}$ over
the children of one particle, the analogue of the moment generating function of a single random
walk; for a single map, $\Lambda(\lambda)=\E\,2^{\lambda X}$. By the theorem of Hammersley, Kingman
and Biggins \cite{Ham74,Kin75,Big76} (see also
\cite{Shi15}), the rightmost particle of such a walk moves at speed
$v=\inf_{\lambda>0}\lambda^{-1}\log_2\Lambda(\lambda)$ bits per generation (here the process
survives forever, since $J(M)$ always contains $r_{S^-}(M)$). The heuristic predicts
that $\pair{S^-}{S^+}$ is unbounded if $v>0$ and bounded if $v<0$; equivalently, when $v\neq 0$,
bounded if and only if $\Lambda(\lambda)<1$ for some $\lambda>0$. Within the model, one direction is elementary:
the expected value of $\sum2^{\lambda\cdot(\text{height})}$ over the particles of generation $t$,
heights being measured from the starting point, is $\Lambda(\lambda)^t$, so by Markov's inequality
the expected number of particles of generation $t\ge 1$ at height at least $0$ is at most
$\Lambda(\lambda)^t$. If $\Lambda(\lambda)<1$, the expected total number of descendants that are at
least as large as the starting value is at most $\Lambda(\lambda)/(1-\Lambda(\lambda))<\infty$; this
is the property checked numerically in \cref{sec:03} below. The heuristic treats the subtrees rooted at different
children as independent, although the children $\Syr(x)$ and $3\Syr(x)$ of $x$ are determined by
the same binary digits of $x$.

\begin{table}[t]
\centering\small
\caption{Branching-random-walk heuristic for families containing the Collatz map:
expected number of children per odd step, $\Lambda_{\min}=\inf_{\lambda>0}\Lambda(\lambda)$, and the
speed $v$ of the rightmost branch (bits per odd step).}
\label{tab:brw}
\begin{tabular}{lcccc}
\toprule
family & $\E\,|J(M)|$ & $\Lambda_{\min}$ & $v$ & prediction\\
\midrule
$\pair{\emptyset}{\{1\}}$ & $5/4$ & $1.2468$ & $+0.639$ & unbounded\\
$\pair{\emptyset}{\{2,3\}}$ & $5/4$ & $1.2044$ & $+0.342$ & unbounded\\
$\pair{\emptyset}{\even}$ & $7/6$ & $1.0986$ & $+0.201$ & unbounded\\
$\pair{\emptyset}{\{2\}}$ & $9/8$ & $1.0778$ & $+0.184$ & unbounded\\
$\pair{\emptyset}{\{3\}}$ & $17/16$ & $0.9942$ & $-0.017$ & bounded\\
$\pair{\emptyset}{\{4\}}$ & $33/32$ & $0.9638$ & $-0.136$ & bounded\\
$\pair{\emptyset}{\{5\}}$ & $65/64$ & $0.9529$ & $-0.214$ & bounded\\
\bottomrule
\end{tabular}
\end{table}

\begin{figure}[t]
\centering
\includegraphics[width=\textwidth]{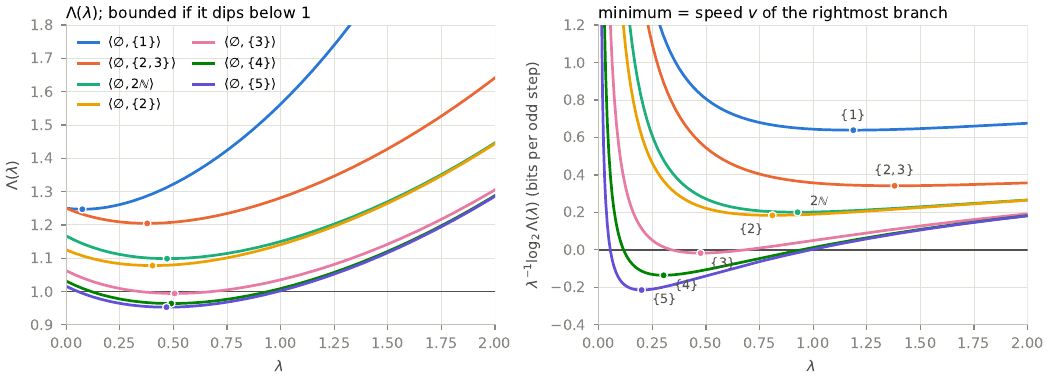}
\caption{The branching-random-walk heuristic for the families of \cref{tab:brw}. Left:
$\Lambda(\lambda)$, with its minimum $\Lambda_{\min}$ marked; a family is predicted to be bounded if
$\Lambda(\lambda)<1$ for some $\lambda>0$. Right: $\lambda^{-1}\log_2\Lambda(\lambda)$, whose minimum is
the speed $v$ of the rightmost branch, labelled by $S^+$. The curves on the left start at
$\Lambda(0)=\E\,|J(M)|>1$, so those on the right tend to $+\infty$ as $\lambda\to 0$.}
\label{fig:brw}
\end{figure}

\Cref{tab:brw} shows the prediction for some families $\pair{\emptyset}{S}$. Every member of these
families is contracting, yet the first four families are predicted to be unbounded. The family
$\pair{\emptyset}{\{3\}}$ is predicted to be bounded, but only barely; by
\cref{cor:collatz} its boundedness would imply that no Collatz orbit diverges. \Cref{fig:brw} shows
the two functions behind the table: $\Lambda(\lambda)$, which dips below $1$ exactly for the three
families predicted to be bounded, and $\lambda^{-1}\log_2\Lambda(\lambda)$, whose minimum is $v$.

\subsection{The family \texorpdfstring{$\pair{\emptyset}{\{1\}}$}{<empty,{1}>}}\label{sec:01}
Here the choice occurs exactly at the odd numbers $x\equiv 1\pmod 8$ (where $M(x)=1$), whose two
successors are $(3x+1)/4$ and $3(3x+1)/4$; every other odd $x$ has the single successor $\Syr(x)$.
Both members, the Collatz map and $T_{\{1\}}$, are contracting.

\emph{Numerics.} Among $n\le 2000$, an exhaustive search shows that $417$ have a finite reachable
set, and for each of the other $1583$ it finds a branch exceeding $2^{300}$ (so these reachable sets
are presumably infinite). Among odd $x<2^{16}$ there are $1345$ with a
finite reachable set (for each of the others the search found a branch exceeding $2^{400}$ or more
than $3\cdot10^5$ reachable numbers). Their number in $[2^{k-1},2^{k})$ grows by a factor of about
$1.5$ per step in $k$, so they appear to have density zero, with roughly $N^{0.58}$ of them below
$N$. A beam
search that keeps the $10^4$ largest numbers in each generation reaches $1861$ bits after $3000$ odd
steps from $27$, that is $0.62$ bits per odd step, against $v=0.639$.

\emph{A positive-drift rule.} We now show rigorously that branches of positive drift exist for
almost every 2-adic starting point. For $D\ge 3$ and a set $\Delta$ of residues modulo $2^D$ that
are $\equiv 1\pmod 8$, let
\[
\psi_\Delta(x)=\begin{cases}
(3x+1)/4, & x\equiv 1\pmod 8 \text{ and } x\bmod 2^D\in\Delta,\\
3(3x+1)/4, & x\equiv 1\pmod 8 \text{ and } x\bmod 2^D\notin\Delta,\\
\Syr(x), & \text{otherwise},
\end{cases}
\]
so that $\psi_\Delta(x)=3^{e(x)-1}(3x+1)/2^{b(x)}$ with $b(x)=v_2(3x+1)$ and $e(x)\in\{1,2\}$. Every
$\psi_\Delta$-orbit of an odd integer consists of odd numbers on one branch of $\pair{\emptyset}{\{1\}}$.

The rule $\psi_\Delta$ depends only on $x\bmod 2^D$. The
next lemma says that, for a random 2-adic $x$, the residue of the next odd number is determined by
the current residue together with fresh, independent binary digits of $x$, so that the residues
form a Markov chain.

\begin{lemma}\label{lem:markov}
Let $x$ be $\mu$-distributed, $x_t=\psi_\Delta^{\,t}(x)$, $s_t=x_t\bmod 2^D$ and $b_t=b(x_t)$.
Then $(s_t)$ is a Markov chain on the $2^{D-1}$ odd residues modulo $2^D$. If $v=v_2(3s+1)<D$,
then from $s$ it moves with probability $2^{-v}$ to each of the $2^{v}$ residues $s'$ with
$s'\equiv 3^{e(s)-1}(3s+1)/2^{v}\pmod{2^{D-v}}$, and $b_t=v$. If $3s+1\equiv 0\pmod{2^D}$, then
$s_{t+1}$ is uniformly distributed and independent of $b_t$, and $\P(b_t=D+i)=2^{-(i+1)}$ for
$i\ge 0$.
\end{lemma}

\begin{proof}
As in \cref{lem:bernoulli}, by induction on $t$: for every finite history
$h=(s_0,b_0,\dots,s_{t-1},b_{t-1},s_t)$ the set of $x$ realizing $h$ is a residue class modulo
$2^{D+b_0+\dots+b_{t-1}}$ that $\psi_\Delta^{\,t}$ maps affinely onto the class $s_t\bmod 2^D$,
multiplying $\mu$ by $2^{b_0+\dots+b_{t-1}}$. Hence $\P(h)=2^{-(D-1)}2^{-(b_0+\dots+b_{t-1})}$, and
the conditional probability of each continuation $(b_t,s_{t+1})$ compatible with $s_t$ is
$2^{-b_t}$, which depends only on $s_t$.
\end{proof}

\begin{theorem}\label{thm:drift}
Let $\psi=\psi_\Delta$ with $D=5$ and $\Delta=\{9\}$, that is: at $x\equiv 1\pmod 8$ take
$(3x+1)/4$ if $x\equiv 9\pmod{32}$ and $3(3x+1)/4$ otherwise. For $\mu$-almost every $x$,
\[
\frac1t\sum_{i<t}e(x_i)\to\frac{20}{17},\qquad \frac1t\sum_{i<t}b(x_i)\to\frac{30}{17},
\]
so that $\frac1t\sum_{i<t}\bigl(e(x_i)\log_2 3-b(x_i)\bigr)\to\frac{10}{17}\log_2\frac98=0.09995\ldots>0$.
\end{theorem}

\begin{proof}
By \cref{lem:markov}, $(s_t)$ is a Markov chain on the $16$ odd residues modulo $32$, with the
transitions listed in \cref{tab:chain}. It is irreducible (from every state one reaches $5$, and
from $5$ every state) and aperiodic ($5\to 5$). (To read the table: from $s=3$ we have
$3s+1=10$, so $b=1$ and the next residue is congruent to $5$ modulo $16$, that is, $5$ or $21$, each
with probability $1/2$; from $s=1$, where $e=2$, the next residue is congruent to $3(3s+1)/4=3$
modulo $8$.) One checks directly that its stationary
distribution is $\pi(s)=8/136$, $9/136$, $6/136$, $11/136$ for $s\equiv 1,3,5,7\pmod 8$
respectively. The process $(s_t,b_t)$ is a positive recurrent Markov chain on a countable set
(by \cref{lem:markov}, the conditional law of $(b_t,s_{t+1})$ given the past depends only on
$s_t$), so by the ergodic theorem the averages converge almost surely to
\begin{align*}
\sum_s\pi(s)\,e(s)&=1+3\cdot\tfrac{8}{136}=\tfrac{20}{17}\qquad\text{(the residues $1,17,25$ have $e=2$)},\\
\sum_s\pi(s)\,\E(b_t\mid s_t=s)&=\tfrac{80}{136}\cdot 1+\tfrac{32}{136}\cdot 2+\tfrac{12}{136}\cdot 3
+\tfrac{6}{136}\cdot 4+\tfrac{6}{136}\cdot 6=\tfrac{240}{136}=\tfrac{30}{17}.
\end{align*}
Finally $\frac{20}{17}\log_2 3-\frac{30}{17}=\frac{10}{17}(2\log_2 3-3)=\frac{10}{17}\log_2\frac98$.
\end{proof}

\begin{figure}[t]
\centering
\includegraphics[width=\textwidth]{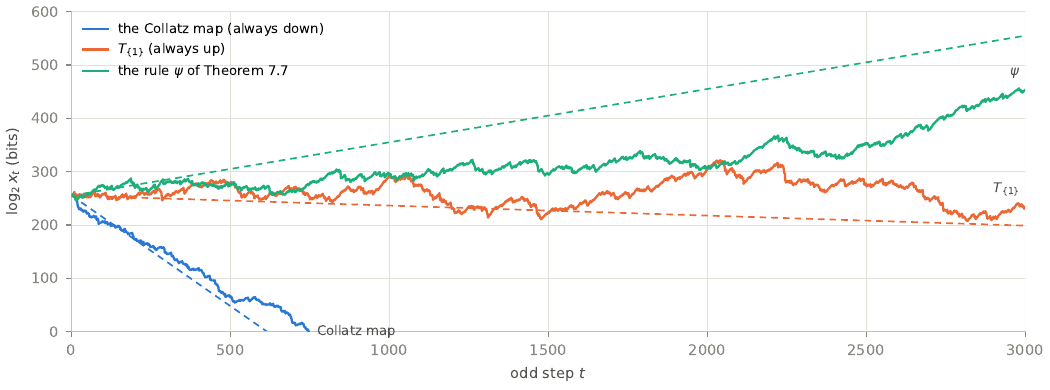}
\caption{Three branches of $\pair{\emptyset}{\{1\}}$ from the same random odd $256$-bit starting value:
the Collatz map (always the step down), $T_{\{1\}}$ (always the step up) and the rule $\psi$ of
\cref{thm:drift}. Dashed: lines with slopes $\log_2 3-2$, $\frac54\log_2 3-2$ and
$\frac{10}{17}\log_2\frac98$. The Collatz orbit reaches $1$ after $746$ odd steps; after $3000$ odd steps
$T_{\{1\}}$ has lost $22$ bits and $\psi$ has gained $198$, against $56$ and $300$ predicted (over
$300$ random $256$-bit starting values the gain of $\psi$ has mean $297$ and standard deviation $82$
bits).}
\label{fig:rule}
\end{figure}

The rule is easy to understand from \cref{tab:chain}. At $x\equiv 9\pmod{32}$ the step up would lead
to a number $\equiv 5\pmod 8$, whose next odd step removes at least three factors $2$ ($b\ge 3$),
while the step down leads to a number $\equiv 7\pmod 8$, which loses only one factor $2$. At
$x\equiv 17\pmod{32}$ the situation is reversed: the step down leads to a number $\equiv 5\pmod 8$
and the step up to a number $\equiv 7\pmod 8$. So the rule uses the choice to steer the orbit away
from residues with large valuations, at the price of giving up a factor $3$ at $x\equiv 9\pmod{32}$.
For example, from $41\equiv 9\pmod{32}$ the step down leads to $31$ and then to $47$, the step up to
$93$ and then to $35$. \Cref{fig:rule} follows the rule and the two members of the family from one
starting value.

\begin{table}[t]
\centering\small
\caption{The Markov chain of \cref{thm:drift}: state $s$ (odd residue modulo $32$), exponent $e$,
valuation $b=v_2(3s+1)$, the states reached (each with probability $2^{-b}$; ``all'' means all $16$
states with probability $1/16$ each), and $136\,\pi(s)$.}
\label{tab:chain}
\begin{tabular}{rccll}
\toprule
$s$ & $e$ & $b$ & next states & $136\,\pi(s)$\\
\midrule
$1$ & $2$ & $2$ & $3,11,19,27$ & $8$\\
$9$ & $1$ & $2$ & $7,15,23,31$ & $8$\\
$17$ & $2$ & $2$ & $7,15,23,31$ & $8$\\
$25$ & $2$ & $2$ & $1,9,17,25$ & $8$\\
$3,\ 11,\ 19,\ 27$ & $1$ & $1$ & $\{5,21\},\ \{1,17\},\ \{13,29\},\ \{9,25\}$ & $9$\\
$7,\ 15,\ 23,\ 31$ & $1$ & $1$ & $\{11,27\},\ \{7,23\},\ \{3,19\},\ \{15,31\}$ & $11$\\
$13$ & $1$ & $3$ & $1,5,9,\dots,29$ (all $\equiv 1\bmod 4$) & $6$\\
$29$ & $1$ & $3$ & $3,7,11,\dots,31$ (all $\equiv 3\bmod 4$) & $6$\\
$5$ & $1$ & $4$ & all & $6$\\
$21$ & $1$ & $\ge 5$ (mean $6$) & all & $6$\\
\bottomrule
\end{tabular}
\end{table}

\begin{remark}\label{rem:lookahead}
The rule of \cref{thm:drift} was found by solving the average-reward Markov decision problem in
which the choice may depend on $x\bmod 2^D$. The optimal drifts are $-0.0188$ for $D=4$ (always
choosing $3(3x+1)/4$, i.e.\ the map $T_{\{1\}}$, is optimal), and
$+0.1000$, $+0.1784$, $+0.1901$, $+0.2156$, $+0.2613$, $+0.2928$, $+0.3185$, $+0.3414$, $+0.3602$,
$+0.3773$ for $D=5,\dots,14$. Thus five bits of the current number suffice for positive drift, and
four do not. For $D=6$ the optimal rule takes $(3x+1)/4$ exactly when $x\bmod 64\in\{1,9,41\}$;
the same proof, with a chain on $32$ states, gives the stationary averages $\frac{20}{17}$ and
$\frac{86}{51}$ and the drift $(60\log_2 3-86)/51=0.1784\ldots$.
\end{remark}

\begin{remark}
Taking the step up whenever some member of the family allows it gives the single map
$T_{S^+}$, and taking the step down whenever possible gives $T_{S^-}$ (\cref{lem:pairs}).
Neither decides whether the family is bounded: for $\pair{\emptyset}{\{1\}}$ both
$T_{\{1\}}$ and the Collatz map are contracting, while the rule of \cref{thm:drift}, which
sometimes takes the step down, has positive drift. A step up changes which steps follow,
so the greedy choice is not the one that grows fastest.
\end{remark}

\begin{theorem}\label{thm:large}
For every $\varepsilon>0$, the set of odd $n$ for which $R_{\pair{\emptyset}{\{1\}}}(n)$ contains
a number larger than $n^{(2/3)\log_2 3-\varepsilon}$ has density one; here
$(2/3)\log_2 3=\log_8 9=1.0566\ldots$ With the $6$-bit rule of \cref{rem:lookahead} the exponent improves to
$(30/43)\log_2 3=1.1057\ldots$; a numerical computation of the stationary averages of the optimal
$14$-bit rule gives about $1.233$.
\end{theorem}

The exponent has
the following origin. The residue of $n$ modulo $2^K$ determines about $K/\bar b$ steps of the
rule, where $\bar b=30/17$ is the average number of factors $2$ removed per step; during these steps
the factors $2$ removed use up the roughly $\log_2 n$ bits of $n$, while the factors $3$ gained
contribute about $\bar e\log_2 3$ bits per step, $\bar e=20/17$. With $K\approx\log_2 n$ this gives
about $(\bar e/\bar b)\log_2 3\cdot\log_2 n$ bits.

\begin{proof}
Let $\bar e=20/17$ and $\bar b=30/17$ be the limits in \cref{thm:drift}, fix $\varepsilon>0$ and
consider odd $n\in(N/2,N]$. Let $t=\lfloor(1-2\varepsilon)\log_2 N/(\bar b+\varepsilon)\rfloor$,
and let $G$ be the set of $x\in\Zu$ with
$\bigl|\sum_{i<t}e(x_i)-t\bar e\bigr|\le\varepsilon t$ and $\bigl|\sum_{i<t}b(x_i)-t\bar b\bigr|\le\varepsilon t$.
By \cref{thm:drift}, $\mu(G)\to 1$. By the proof of \cref{lem:markov}, $G$ is a union of residue
classes modulo $2^{D+\sum_{i<t}b(x_i)}$, hence modulo $2^{D+(1-2\varepsilon)\log_2 N}$, and as in the
proof of \cref{thm:density} the proportion of odd $n\in(N/2,N]$ lying in $G$ tends to $1$. For such
$n$, by \eqref{eq:ratio} applied to each step,
\[
x_t\ \ge\ n\,\frac{3^{\sum e(x_i)}}{2^{\sum b(x_i)}}\ \ge\ 3^{\,t(\bar e-\varepsilon)}
\]
for $N$ large, because $2^{\sum b(x_i)}\le 2^{(1-2\varepsilon)\log_2 N}<n$. Since
$t\ge(1-3\varepsilon)\log_2 n/(\bar b+\varepsilon)$ for large $N$, we get
$\log_2 x_t\ge\frac{(1-3\varepsilon)(\bar e-\varepsilon)}{\bar b+\varepsilon}\log_2 3\cdot\log_2 n$, and
$\bar e/\bar b=2/3$. The number $x_t$ lies on a branch from $n$. The other exponents are obtained in
the same way from the stationary averages of the corresponding rules.
\end{proof}

\Cref{thm:large} does not show that $\pair{\emptyset}{\{1\}}$ is unbounded: the large numbers it
produces might all belong to the exceptional set of density zero. The following computation shows
that residue classes alone cannot settle the question.

\begin{proposition}\label{prop:rational-traps}
Every odd residue class modulo $2^{10}$, and hence every odd residue class modulo $2^L$ for
$L\le 10$, contains a positive rational number
$y=p/q$, $q$ odd, whose reachable set under $\pair{\emptyset}{\{1\}}$ is finite. The denominators
can be chosen at most $2187$ and the reachable sets have at most $73$ odd elements. For example,
$41/9\equiv 33\pmod{256}$, and
\[
\tfrac{41}{9}\ \mapsto\ \tfrac{11}{3}\ \text{or}\ 11,\qquad \tfrac{11}{3}\ \mapsto\ 3\ \text{or}\ 9,
\]
where $3$, $9$ and $11$ have finite reachable sets; $R(41/9)$ has $24$ odd elements.
Explicitly, the odd elements of $R(41/9)$ are $41/9$, $11/3$ and the $22$ integers $1$, $3$, $5$,
$7$, $9$, $11$, $13$, $17$, $19$, $21$, $29$, $39$, $59$, $67$, $85$, $89$, $101$, $151$, $201$,
$227$, $341$, $453$. Thus $41/9$ lies in the residue class of $33$ modulo $2^8$ but has a finite
reachable set, while integers with a finite reachable set appear to be rare.
\end{proposition}

\begin{proof}
For a rational $y=p/q$ the odd-to-odd steps of the family act on the numerator as
$p\mapsto(3p+q)/2^{v_2(3p+q)}$, with the extra choice $p\mapsto 3(3p+q)/4$ when
$v_2(3p+q)=2$. For each class we found $p,q$ by a computer search and then enumerated the reachable
set with exact rational arithmetic. The lists of witnesses are in the ancillary files.
\end{proof}

To state the consequence, call a nonempty set $C$ of odd residue classes modulo $2^L$ a
\emph{growth certificate} if there are $\gamma>1$ and, for every $c\in C$, a finite decision
tree that reads further binary digits of $x\in c$ and chooses moves of the family, such that
every $x\in c$ is led along a branch to a number $x'$ with $x'\bmod 2^L\in C$ and
$x'\ge\gamma x$, where the last inequality is derived from the multipliers $3^{e}/2^{b}$ and
\eqref{eq:ratio}, so that it holds for all positive $x$ in $c$. A growth certificate would show
that every positive integer in a class of $C$ has an unbounded branch. In other words, a
growth certificate is a finite strategy, based on residues, that guarantees growth by the factor
$\gamma$ from every number in a class of $C$ to a number in a class of $C$.

\begin{corollary}\label{cor:nocert}
There is no growth certificate for $\pair{\emptyset}{\{1\}}$ modulo $2^L$ with $L\le 10$.
\end{corollary}

\begin{proof}
A certificate modulo $2^L$ with $L\le 10$ is also one modulo $2^{10}$. Applied to the rational $y$ of
\cref{prop:rational-traps} in a class of $C$, it would produce numbers $\ge\gamma^k y$ in $R(y)$
for every $k$, contradicting finiteness.
\end{proof}

For $L=12$ we found such rational witnesses in $501$ of the $512$ classes $\equiv 1\pmod 8$. A
direct search for certificates in which the trees read at most $8$ extra digits and make at most
$40$ odd steps, landing on classes $\equiv 1\pmod 8$, found none for $L=12$ and $L=14$: the largest
self-sustaining set of classes is empty. Within these bounds the search is exhaustive; in particular
it never reached its limit on the number of nodes ($2\cdot10^6$ per class).

\subsection{The family \texorpdfstring{$\pair{\emptyset}{\{3\}}$}{<empty,{3}>}}\label{sec:03}
Here the choice occurs at odd $x$ with $M(x)=3$, i.e.\ $x\equiv 5\pmod{32}$, between $\Syr(x)$ and
$3\Syr(x)$. We checked by induction that every $n<10^8$ has a finite reachable set: for every odd
$n$, the set of numbers $\ge n$ reachable from $n$ without passing below $n$ is finite. This
suffices: by induction on $n$, every branch from $n$ either stays in this finite set or passes below
$n$, and every branch from an even number passes below it within two steps. The largest
such set has $2431$ elements (for $n=37{,}537{,}407$), and the largest climb is $31$ bits (from
$n=4{,}196{,}607$). In a simulation of the branching random walk of \cref{sec:brw}, the probability
that the maximum climbs at least $30$ bits before all branches fall below the start is about
$10^{-7}$, so for $5\cdot 10^7$ odd starting values one expects a largest climb of $30$ to $34$
bits, in agreement with the data.

\section{Unified maps}\label{sec:unified}

The maps $T_S$ and $T^-_S$ agree at even numbers and differ at odd numbers, which all go up under
$T_S$ and all go down under $T^-_S$. We now let the direction at odd numbers depend on a 2-adic
valuation as well. For odd $x$ we use $v_2(x+1)$, the number of binary digits $1$ at the end of
$x$; it is the exponent $M^-(x)+1$ of \cref{sec:minus}. The question is whether the maps between the
two families behave in a new way. The computations below suggest that they do not, and that a single
exponent decides: if the odd numbers $x\equiv 1\pmod 4$, those with $v_2(x+1)=1$, go down, a map
behaves like the maps $T^-_S$ (\cref{sec:unified-down}), and if they go up, like the maps $T_S$
(\cref{sec:unified-up}).

\begin{definition}\label{def:unified}
For $\So,\Se\subseteq\N$ let $T_{\So,\Se}\colon\N\to\N$ be
\begin{equation}\label{eq:unified}
T_{\So,\Se}(x)=\begin{cases}
(3x+1)/2, & x \text{ odd},\ m\in\So,\\
(x+1)/2, & x \text{ odd},\ m\notin\So,\\
3x/2, & x \text{ even},\ m\in\Se,\\
k, & x \text{ even},\ m\notin\Se,
\end{cases}
\end{equation}
where $x+1=2^mk$ if $x$ is odd and $x=2^mk$ if $x$ is even, with $m\ge 1$ and $k$ odd.
\end{definition}

Thus $\So$ is the set of exponents at which odd numbers go up and $\Se$ the set of exponents at
which even numbers go up, and
\[
T_S=T_{\N,S},\qquad T^-_S=T_{\emptyset,S}.
\]
We call $T_{\So,\Se}$ \emph{mixed} if $\So\neq\emptyset$ and $\So\neq\N$; the mixed maps belong to
neither family. Since $(3x+1)/2>(x+1)/2$ and $3x/2>k$, \cref{lem:monotone} extends:
$T_{\So,\Se}\le T_{\So',\Se'}$ pointwise if $\So\subseteq\So'$ and $\Se\subseteq\Se'$. The
family runs from $T_{\emptyset,\emptyset}=T^-_\emptyset$, whose orbits all decrease to $1$, to
$T_{\N,\N}=T_\N$, whose orbits all diverge, and it contains the Collatz map twice: as
$T_{\N,\emptyset}=T_\emptyset$ and, after the shift by one of \cref{prop:shift}, as
$T_{\emptyset,\N}=T^-_\N$. For example,
\begin{align*}
T_{\{2\},\emptyset}&:\quad 27,\,41,\,21,\,11,\,17,\,9,\,5,\,3,\,5,\dots\\
T_{\{1\},[1,3]}&:\quad 37,\,56,\,84,\,126,\,189,\,284,\,426,\,639,\,320,\,5,\,8,\,12,\,18,\,27,\,14,\,21,\,32,\,1,\,2,\,3,\,2,\dots
\end{align*}
The first orbit ends in the cycle $(3,5)$, in which $3$ goes up ($v_2(4)=2\in\So$) and $5$ goes
down ($v_2(6)=1\notin\So$); so $(3,5)$ is a cycle of no map $T_S$ and of no map $T^-_S$. The
second orbit ends in the cycle $(2,3)$ of \cref{sec:minus}. \Cref{rem:switch} below explains why
the two cycles are closely related.

\subsection{The block map}\label{sec:unified-block}
Since $(3x+1)/2+1=\frac32(x+1)$, a step up acts in the same way at odd and at even numbers: in the
notation of \cref{def:unified} it replaces $2^mk$ by $2^{m-1}\cdot 3k$, and the parity of the
number is kept as long as $m\ge 2$ (for $m=1$ an odd $x$ goes up to the even number $3k-1$, and an
even $x$ to the odd number $3k$). So at odd numbers, too, a run of steps up passes through the
exponents $m,m-1,\dots$, and the run of consecutive elements of $\So$ that ends at $m$ decides how
many steps up are taken. Under $T_S$ every odd $x$ goes up, and $x$ is followed by $m-1$ further odd
numbers. A step down at an odd $x$ gives $(x+1)/2=2^{m-1}k$: it lowers the exponent by one, as a step
up does, and for $m\ge 2$ hands the even number $2^{m-1}k$ over to the rules for even numbers. A
step down at an even number gives its odd part at once.

\begin{lemma}\label{lem:Uunified}
Let $x$ be odd, $x+1=2^mk$ with $k$ odd, and $r=r_{\So}(m)$. The $T_{\So,\Se}$-orbit of $x$
begins with the $r$ steps up $x_i\mapsto x_{i+1}$, $0\le i<r$, where $x_i=2^{m-i}3^ik-1$. Let
\[
y=\begin{cases}(x_r+1)/2=2^{m-r-1}3^rk, & r<m,\\ x_r=3^mk-1, & r=m,\end{cases}
\]
which is the next term of the orbit if $r<m$ and equal to $x_r$ if $r=m$, and write $y=2^Mk'$ with
$M\ge 0$ and $k'$ odd. Then the orbit of $x$ reaches the odd number
\begin{equation}\label{eq:Uunified}
U_{\So,\Se}(x)=3^{\,r_{\Se}(M)}\,k'
\end{equation}
after $\ell=r+[\,r<m\,]+r_{\Se}(M)+[\,r_{\Se}(M)<M\,]$ steps; all terms after $x_r$ and before
$U_{\So,\Se}(x)$ are even. Explicitly,
\[
U_{\So,\Se}(x)=\begin{cases}
3^{\,r+r_{\Se}(m-r-1)}\,(x+1)/2^m, & r<m,\\[2pt]
3^{\,r_{\Se}(M)}\,\bigl(3^m(x+1)/2^m-1\bigr)/2^M, & r=m.
\end{cases}
\]
\end{lemma}

\begin{proof}
By the rule above, the step up at $x_i$ is taken if $m-i\in\So$, and it leads to $x_{i+1}$. By
\eqref{eq:r}, exactly $r$ steps up are taken before the exponent reaches $m-r$, which is $0$ if
$r=m$ and not in $\So$ if $r<m$; in the second case $x_r$ is odd and goes down to $y$. If $M=0$,
then $y=k'$ is odd and $\ell=r+1$. If $M\ge1$, then by the last sentence of \cref{lem:US} the even
number $y=2^Mk'$ reaches the odd number $3^{r_{\Se}(M)}k'$ after $r_{\Se}(M)+[\,r_{\Se}(M)<M\,]$
steps, through even numbers.
\end{proof}

We call $U_{\So,\Se}$ the \emph{block map}, and the part of an orbit from an odd number $x$ to
$U_{\So,\Se}(x)$ a \emph{block}. For $\So=\emptyset$ we have $r=0$, $y=(x+1)/2$ and $M=M^-(x)$, so
$U_{\emptyset,S}=U^-_S$ by \eqref{eq:Uminus}. For $\So=\N$ we have $r=m$, and $x_0,\dots,x_{m-1}$
are consecutive odd numbers of the $T_S$-orbit, with $M(x_i)=0$ for $i<m-1$ and
$(3x_{m-1}+1)/2=y=2^{M}k'$, so that $M(x_{m-1})=M$ and $U_{\N,S}(x)=U_S^{\,m}(x)$: a block of $T_S$
consists of $v_2(x+1)$ odd steps. In general $y=2^Mk'$ is the number at which the rules for even
numbers take over; in the two families it is $(3x_{m-1}+1)/2=2^M\Syr(x_{m-1})$, as in
\cref{sec:maps}, and $(x+1)/2$, as in \cref{sec:minus}. As in \cref{sec:rationals},
\eqref{eq:unified} and \eqref{eq:Uunified} also make sense for positive rational numbers with odd
denominator, and \eqref{eq:Uunified} for all $x\in\Zu$ outside a countable set.

\begin{remark}[Switching]\label{rem:switch}
In the case $r<m$, a block consists of $r$ steps up at odd numbers, one step down at an odd number,
and the steps at even numbers that follow: $b=r_{\Se}(M)$ steps up and, if $b<M$, a step down to the
odd part. Every step up multiplies $x_i+1$, or the even number, by $\frac32$, and the step down at
the odd number $x_r$ halves $x_r+1$. So the block multiplies $x+1$ by $3^a/2^m$, where $a=r+b$ is
the number of steps up: for a given $a$, it does not matter where the step down at an odd number is
taken. Moving this step within the block, by trading steps up at odd numbers for steps up at even
numbers or conversely, changes neither the end of the block nor its number of steps, but only the
numbers passed through, which are of the form $2^i3^hk-1$ before the step down and $2^i3^hk$ after
it. For example, $x=3$ has $m=2$, and its block with $a=1$ is $3\to 5\to 3$ if $2\in\So$ and
$1\notin\So$, and $3\to 2\to 3$ if $2\notin\So$ and $1\in\Se$: these are the cycles $(3,5)$ and
$(2,3)$. Consequently, different pairs $(\So,\Se)$ can have the same block map; $(\{1\},[1,3])$
and $(\{1,3\},[1,3])$ are an example. We say that two cycles are in the same \emph{class} if they
differ only by such moves.
\end{remark}

\subsection{Cycles and a census}\label{sec:unified-census}
A cycle $C$ of a map $T_{\So,\Se}$ determines four finite sets of exponents: the sets
$U_{\mathrm o}(C)$ and $D_{\mathrm o}(C)$ of the exponents $v_2(x+1)$ of the odd elements $x\in C$
that go up and down, and the sets $U_{\mathrm e}(C)$ and $D_{\mathrm e}(C)$ of the exponents
$v_2(x)$ of the even elements that go up and down. As in \cref{prop:cycle-interval}, $C$ is a cycle
of $T_{\So,\Se}$ if and only if
\[
U_{\mathrm o}(C)\subseteq\So\subseteq\N\setminus D_{\mathrm o}(C)\quad\text{and}\quad
U_{\mathrm e}(C)\subseteq\Se\subseteq\N\setminus D_{\mathrm e}(C).
\]
The $+$ cycles of \cref{app:cycles} are the cycles with $D_{\mathrm o}(C)=\emptyset$, with
$U(C)=U_{\mathrm e}(C)$ and $D(C)=D_{\mathrm e}(C)$, and the $-$ cycles are those with
$U_{\mathrm o}(C)=\emptyset$. Both kinds also occur in mixed maps: $(1,2)$, for instance, is a cycle
of $T_{\So,\Se}$ whenever $1\in\So$ and $1\notin\Se$, and $(1)$ whenever $1\notin\So$. We call $C$
\emph{mixed} if $U_{\mathrm o}(C)$ and $D_{\mathrm o}(C)$ are both nonempty; mixed cycles occur only
in mixed maps. The simplest one is $(3,5)$, with $U_{\mathrm o}=\{2\}$, $D_{\mathrm o}=\{1\}$ and no
even elements, a cycle of $T_{\So,\Se}$ whenever $2\in\So$ and $1\notin\So$; by \cref{rem:switch} it
is in the class of the $-$ cycle $(2,3)$.

To find the cycles of all unified maps at once we extend the search of \cref{app:cycles}. For each
$n$ we follow the orbits from $n$ in which every exponent $v_2(x+1)$ of an odd number and every
exponent $v_2(x)$ of an even number is assigned the step up or the step down when it first occurs,
following both assignments; a branch is abandoned when it drops below $n$, when it exceeds a fixed
bound ($2^{120}$ or $2^{64}$ below), or when it enters a loop that does not contain $n$. This finds,
exactly once, every cycle of every map $T_{\So,\Se}$ whose smallest element is $n$ and whose largest
element is below the bound. We searched the pairs with $1\notin\So$ and those with $1\in\So$
separately; a cycle that does not meet the exponent $1$ at an odd number, such as $(2,3)$, is found in
both searches. As a check, the $+$ cycles found are exactly those of \cref{tab:census} below the same
bounds, the $-$ cycles are $(1)$ and $(2,3)$, and every cycle was verified again, with its four sets,
by an independent program. For counting, we also group the cycles into classes: two cycles are in the
same class if and only if they start their blocks at the same odd numbers, in the same order and with
the same case of \cref{lem:Uunified} in each block.

\subsection{The random model}\label{sec:unified-model}
The blocks of a random 2-adic integer are independent, as the odd steps are in
\cref{lem:bernoulli}. The odd $x$ with $v_2(x+1)=m$ form the residue class of $2^m-1$ modulo
$2^{m+1}$, of measure $2^{-m}$. If $r<m$, the block map $x\mapsto 3^a(x+1)/2^m$ maps this class
affinely onto $\Zu$. If $r=m$, then $k=(x+1)/2^m$ is uniformly distributed in $\Zu$ and $y=3^mk-1$ in
$2\Zt$, so the exponent $M=v_2(y)$ is new, with $\P(M=i)=2^{-i}$ for $i\ge1$; each set
$\{v_2(x+1)=m,\ v_2(y)=M\}$ is a residue class modulo $2^{m+M+1}$, which $U_{\So,\Se}$ maps affinely
onto $\Zu$. The induction in the proof of \cref{lem:bernoulli} then shows that for $\mu$-distributed
$x$ the exponents $m_i=v_2(x_i+1)$ along the orbit $x_{i+1}=U_{\So,\Se}(x_i)$ are independent, with
$\P(m_i=m)=2^{-m}$, and that the exponents $M$ of the blocks with $r=m$ are independent of them and of
each other, with $\P(M=i)=2^{-i}$.

A block takes $u=r+r_{\Se}(M)$ steps up and removes $h$ factors $2$, where $h=m$ if $r<m$ and
$h=m+M$ if $r=m$; up to a correction of order $1/x$, it changes $\log_2 x$ by $u\log_2 3-h$. Since
the blocks have different lengths, we count time in halvings, each of which uses up one binary digit
of $x$ (the set of $x$ with a given block has measure $2^{-h}$), and measure the drift in bits per two
halvings, the average number of halvings in an odd step of $T_S$ and of $T^-_S$:
\begin{equation}\label{eq:drift-unified}
\delta(\So,\Se)=2\bigl(\sigma\log_2 3-1\bigr),\qquad \sigma=\sigma(\So,\Se)=\frac{\E u}{\E h}.
\end{equation}
Here $\E h=2+2\pi$, where $\pi=\sum_{m:\,[1,m]\subseteq\So}2^{-m}$ is the probability of the case
$r=m$, and $\E u=\sum_{m\ge1}2^{-m}\bar u(m)$, where $\bar u(m)=r+r_{\Se}(m-r-1)$ if
$r=r_{\So}(m)<m$, and $\bar u(m)=m+2\rho(\Se)$ if $r=m$. For $\So=\N$ we get $\E h=4$ and
$\E u=2+2\rho(\Se)$, and for $\So=\emptyset$ we get $\E h=2$ and $\E u=\rho(\Se)$; so
$\delta(\N,S)=\delta(S)$ and $\delta(\emptyset,S)=\delta^-(S)$. The sign of $\delta$ is the sign of
the drift per block, $\E u\log_2 3-\E h$, so $T_{\So,\Se}$ is contracting if and only if
$\sigma<\log_3 2=0.6309\ldots$, that is, if a block takes on average fewer than $0.6309$ steps up per
halving. For example, $\sigma(\{1\},[1,3])=5/8=\sigma(\N,\{1\})$, so $T_{\{1\},[1,3]}$ has the
drift $-0.0188$ of $T_{\{1\}}$.

\begin{proposition}\label{prop:sigma-monotone}
If $\So\subseteq\So'$ and $\Se\subseteq\Se'$, then $\sigma(\So,\Se)\le\sigma(\So',\Se')$. In
particular, among all maps that contain a given cycle $C$, the map $T_{U_{\mathrm o}(C),U_{\mathrm
e}(C)}$ has the smallest drift.
\end{proposition}

\begin{proof}
It suffices to add one element $s$ to one of the sets. If $s$ is added to $\Se$, then no $\bar u(m)$
decreases and $\E h$ does not change. Let $s$ be added to $\So$, and let $c$ be the smallest
positive integer not in $\So$, so that $s\ge c$. In the blocks with $r=r_{\So}(m)<m$ that remain in
this case, $r$ can only grow, and $\bar u(m)=r+r_{\Se}(m-r-1)$ does not decrease when $r$ grows,
because $r_{\Se}(i)\le r_{\Se}(i-d)+d$ for $0\le d\le i$; if $s>c$, these are the only changes, and
$\E h$ does not change. If $s=c$, the blocks with $c\le m<c'$, where $c'>c$ is the next integer not
in $\So$ ($c'=\infty$ if there is none), pass from the case $r<m$, with $\bar u(m)\le m-1$ and $h=m$,
to the case $r=m$, with $\bar u(m)=m+2\rho(\Se)$ and on average $m+2$ halvings. So $\E u$ and $\E h$ grow by
amounts whose ratio is at least $\frac12+\rho(\Se)$. Before the change, $\bar u(m)\le m-1$ for $m\ge
c$ gives $\E u\le 1+W+2\rho(\Se)W$ with $W=1-2^{1-c}$, while $\E h=2+2W$, so that $\sigma\le
\frac12+\rho(\Se)W/(1+W)\le\frac12+\rho(\Se)$. Adding to the numerator and the denominator of a
fraction amounts whose ratio is at least the fraction does not decrease it.
\end{proof}

The drift per block, on the other hand, is not monotone: it is $-2$ for $T_{\emptyset,\emptyset}$
and $\frac12\log_2 3-3=-2.21$ for $T_{\{1\},\emptyset}$, whose blocks with $m=1$ are longer.

\subsection{Odd numbers \texorpdfstring{$x\equiv 1\pmod 4$}{x = 1 mod 4} going down}\label{sec:unified-down}
If $1\notin\So$, the odd numbers $x\equiv 1\pmod 4$ go down to $(x+1)/2$. Then $[1,m]\not\subseteq\So$
for every $m$, so every block is of the case $r<m$ of \cref{lem:Uunified}, as for the maps $T^-_S$,
and the block map has the form of \eqref{eq:Uminus}, with a more general exponent of $3$.

\begin{proposition}\label{prop:down}
Let $1\notin\So$, and for $m\ge1$ let $a(m)=r_{\So}(m)+r_{\Se}\bigl(m-r_{\So}(m)-1\bigr)$.
\begin{enumerate}[label=\textup{(\alph*)}]
\item For odd $x$ with $v_2(x+1)=m$ we have $U_{\So,\Se}(x)=3^{a(m)}(x+1)/2^m$ and $a(m)\le m-1$.
Hence $U_{\So,\Se}(x)\le U^-_\N(x)$ for every odd $x$, and $\delta(\So,\Se)\le\log_2 3-2$.
\item $a(m)=m-1$ for every $m$ if and only if $[1,g-1]\subseteq\Se$ for every $g\in\N\setminus\So$.
In this case $U_{\So,\Se}=U^-_\N$, so every $T_{\So,\Se}$-orbit is eventually periodic if and only
if no Collatz orbit diverges, and every cycle of $T_{\So,\Se}$ is in the class of a cycle of
$T^-_\N$, that is, of a Collatz cycle shifted by one.
\item If $3^{a(m)}<2^m$ for every $m$, then every $T_{\So,\Se}$-orbit of a positive integer
reaches the fixed point $1$ or a cycle through $3$, which is $(3,5)$ if $2\in\So$, and $(2,3)$ if
$2\notin\So$ and $1\in\Se$.
\end{enumerate}
\end{proposition}

\begin{proof}
(a) Since $1\notin\So$, we have $r=r_{\So}(m)<m$, so \cref{lem:Uunified} gives the formula, and
$r_{\Se}(m-r-1)\le m-r-1$. For $(\So,\Se)=(\emptyset,\N)$ we get $a(m)=m-1$, that is,
$U^-_\N(x)=3^{m-1}(x+1)/2^m$. Finally $\E h=2$ and $\E u=\sum_m2^{-m}a(m)\le\sum_m2^{-m}(m-1)=1$.

(b) Let $g=m-r_{\So}(m)$, the largest element of $[1,m]$ not in $\So$. Then
$a(m)=m-g+r_{\Se}(g-1)$, which equals $m-1$ if and only if $r_{\Se}(g-1)=g-1$, that is,
$[1,g-1]\subseteq\Se$; and every $g\notin\So$ occurs, for $m=g$. If the block maps agree, the orbits
of $T_{\So,\Se}$ and $T^-_\N$ through an odd number start their blocks at the same odd numbers, and
corresponding blocks differ only by the moves of \cref{rem:switch}. A block has finitely many steps,
so an orbit is eventually periodic if and only if the sequence of its block starts is; and by
\cref{prop:shift}, $T^-_\N$ is the Collatz map shifted by one.

(c) Write $x+1=2^mk$. Then $U_{\So,\Se}(x)=3^{a(m)}k<2^mk=x+1$, so $U_{\So,\Se}(x)\le x$, with
equality if and only if $k(2^m-3^{a(m)})=1$, that is, $k=1$ and $2^m-3^{a(m)}=1$. For $m\ge3$ this is
impossible, since $2^m-1\equiv 7\pmod 8$ while $3^{a}\equiv 1$ or $3\pmod 8$. So either $m=1$,
$a(1)=0$ and $x=1$, or $m=2$, $a(2)=1$ and $x=3$. Hence the odd numbers at which blocks start
decrease strictly until they reach $1$ or $3$. The block of $1$ is $1\to1$, and the block of $3$ is
$3\to5\to3$ if $2\in\So$, $3\to2\to3$ if $2\notin\So$ and $1\in\Se$, and $3\to2\to1$ otherwise.
\end{proof}

Part (a) says that without the step up at the exponent $1$, no block multiplies $x+1$ by more than
the corresponding block of the Collatz map $T^-_\N$, which goes down once and then up at every
exponent of the even numbers that follow, the largest possible number $m-1$ of steps up. So every map
with $1\notin\So$ is at least as contracting as the Collatz map; its drift lies between
$\delta(\emptyset,\emptyset)=-2$ and $\log_2 3-2=-0.415$.

Part (b) gives infinitely many maps that are the Collatz map in disguise, in a different way than in
\cref{prop:shift}: for instance $T_{[2,\infty),S}$ for every $S$, and $T_{[3,\infty),S}$ for every
$S$ with $1\in S$. Under $T_{[2,\infty),S}$, an odd $x$ with $x+1=2^mk$, $m\ge2$, goes up through the
odd numbers $2^{m-i}3^ik-1$ and down at the exponent $1$, to $3^{m-1}k$; under $T^-_\N$ it goes down
first, to $2^{m-1}k$, and then up through even numbers, to the same odd number $3^{m-1}k$. The cycle
$(3,5)$ of $T_{[2,\infty),S}$ thus corresponds to the cycle $(2,3)$ of $T^-_\N$.

Part (c) contains \cref{prop:simple}, the case $\So=\emptyset$. It applies, for instance, to the first
example of this section, $T_{\{2\},\emptyset}$, for which $a(2)=1$ and $a(m)=0$ for $m\ne2$: every
orbit ends in $(1)$ or in $(3,5)$. When the criterion fails, $3^{a(m)}>2^m$ for some $m$, and as in
\cref{sec:criterion} the negative rational number $x^*=-3^{a(m)}/(3^{a(m)}-2^m)$, which satisfies
$v_2(x^*+1)=m$, repeats this increasing block forever; so the descent of all orbits cannot be proved
from residue classes alone.

The census of \cref{sec:unified-census}, run over all pairs with $1\notin\So$ and all $n<10^9$,
never reached the bound $2^{120}$ and found only the cycles $(1)$, $(2,3)$ and $(3,5)$. So for every
$\So$ with $1\notin\So$ and every $\Se$, these are the only cycles of $T_{\So,\Se}$ with smallest
element below $10^9$; up to switching, they are the two shifted Collatz cycles. This includes the maps
$T^-_S$, for which \cref{app:cycles} went up to $10^6$. In this respect, too, the maps with
$1\notin\So$, mixed or not, behave like the maps $T^-_S$.

\subsection{Odd numbers \texorpdfstring{$x\equiv 1\pmod 4$}{x = 1 mod 4} going up}\label{sec:unified-up}
If $1\in\So$, let $c\ge2$ be the smallest positive integer not in $\So$ ($c=\infty$ if $\So=\N$).
The blocks with $m<c$ are of the case $r=m$: they go up at all the odd exponents $m,m-1,\dots,1$ and
reach the even number $3^m(x+1)/2^m-1$, whose exponent $M$ is new, as after the step $x\mapsto 3x+1$
of the Collatz map. Such a block multiplies $x+1$ by about $3^{m+r_{\Se}(M)}/2^{m+M}$, which for
$M=1$ is larger than the factor $3^{m-1}/2^m$ of the corresponding block of the Collatz map $T^-_\N$.
So these maps can be expanding: by \cref{prop:sigma-monotone} the drift of the maps with $1\in\So$
ranges from $\delta(\{1\},\emptyset)=\frac13\log_2 3-2=-1.47$ to $\delta(\N,\N)=2\log_2 3-2=1.17$, and
the maps $T_S$ are among them.

Here the census finds many cycles. Below $1000$ (with largest element below $2^{120}$) it finds the
$43$ $+$ cycles of \cref{tab:census}, the cycle $(2,3)$, which does not meet the exponent $1$ at an
odd number, and $17{,}200$ mixed cycles, which form $2{,}436$ classes. Their number grows with the
bound on the smallest element, but slowly, like that of the $+$ cycles: with largest element below
$2^{64}$ there are $217$, $650$ and $845$ classes of mixed cycles with smallest element below
$10^2$, $10^3$ and $10^4$, against $24$, $37$ and $45$ $+$ cycles.

\begin{figure}[!htbp]
\centering
\includegraphics[width=\textwidth]{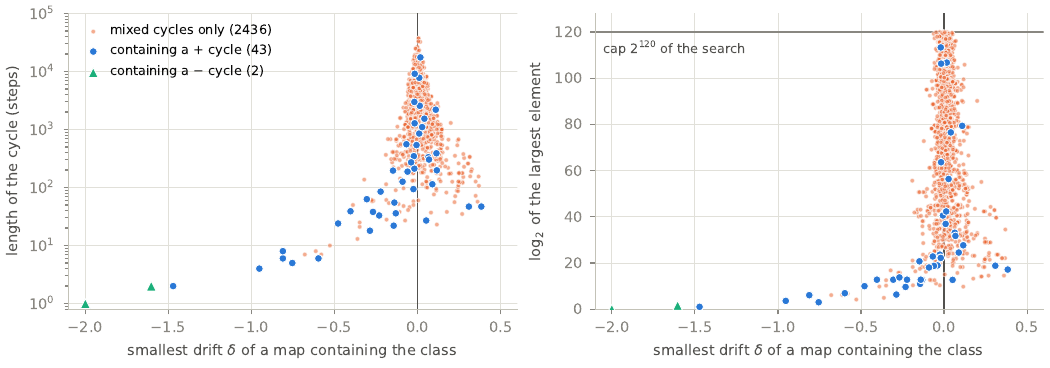}
\caption{The classes of cycles found by the census of \cref{sec:unified-census}: all cycles with
smallest element below $1000$ and largest element below $2^{120}$ of the maps with $1\in\So$, and the
three cycles of the maps with $1\notin\So$. Length (left) and $\log_2$ of the largest element
(right), against the drift of the most contracting map that contains the class
(\cref{prop:sigma-monotone}). Orange: classes of mixed cycles; blue: $+$ cycles; green: the classes
of $(1)$ and of $(2,3)$ and $(3,5)$. The vertical line is the threshold $\delta=0$, the horizontal
line on the right the bound $2^{120}$ of the search.}
\label{fig:unified}
\end{figure}

\Cref{fig:unified} plots the classes against the drift of the most contracting map that contains
them. The mixed classes crowd around the threshold: for $2{,}253$ of the $2{,}436$ ($92\%$) this most
contracting map has $|\delta|<0.1$; these classes have a median length of about $3200$ steps, and
$93\%$ of them have elements above $2^{40}$. This is the picture of \cref{fig:census}, much amplified. Near the threshold orbits
make long excursions (\cref{sec:predictions}), so that long cycles with large elements can close up,
and the unified maps, with two sets as parameters, include far more maps near the threshold than the
family $T_S$.

So the large number of mixed cycles reflects the size of the family rather than the behaviour of
single maps. To compare single maps, we read off from the census the cycles of $512$ maps $T_S$ and
of $32{,}512$ mixed maps with $1\in\So$ and compare them at equal drift (\cref{tab:unified}). A mixed
map has about as many cycles as a map $T_S$ with the same drift, two or three with smallest element
below $1000$ near the threshold. At drifts below $-0.5$, which no map $T_S$ has, the mixed maps have
one to three cycles each, all small (with at most ten elements, the largest at most $120$), but not
only the shifted Collatz cycles: besides $(1,2)$, $(2,3)$ and $(1,2,3,5,8)$ they include, for example,
$C^+_{\{3\}}(3)=(3,5,8,12)$, $C^+_{\{2\}}(19)$ and the mixed cycle $(17,26,39,59,30,45,68)$, which occurs
in maps with drift down to $-0.68$. The maps with $1\notin\So$ have no such cycles at any drift. We
also checked
$22$ mixed maps with drifts between $-1.47$ and $-0.019$ directly, by the method of
\cref{sec:method}, for all $n<10^6$: every orbit is eventually periodic, and each map has one to five
cycles with smallest element below $10^6$, as the maps of \cref{tab:cycles}. For instance, under
$T_{\{1\},[1,3]}$, which has the drift of $T_{\{1\}}$, every $n<10^6$ ends in $(2,3)$, and the
longest glide ($114{,}654$ steps) and the highest peak ($851$ bits) are of the same size as those of
$T_{\{1\}}$ for $n<10^6$ ($86{,}626$ steps and $602$ bits).

\begin{table}[!htbp]
\centering\small
\caption{Cycles of single maps, read off from the census of \cref{sec:unified-census}: for the
contracting maps with drift $\delta$ in the given range, the number of maps and the mean number of
cycles with smallest element below $1000$ and largest element below $2^{120}$. Maps $T_S$: all $S$
for which $S\cap[9,\infty)$ is empty or equal to $[9,\infty)$. Mixed maps: all pairs with
$1\in\So\neq\N$ for which $\So\cap[8,\infty)$ and $\Se\cap[8,\infty)$ are empty or equal to
$[8,\infty)$.}
\label{tab:unified}
\begin{tabular}{lrrrr}
\toprule
 & \multicolumn{2}{c}{maps $T_S$} & \multicolumn{2}{c}{mixed maps with $1\in\So$}\\
\cmidrule(lr){2-3}\cmidrule(lr){4-5}
$\delta$ & maps & cycles per map & maps & cycles per map\\
\midrule
$[-1.5,-1.0)$ & -- & -- & $4{,}590$ & $1.00$\\
$[-1.0,-0.5)$ & -- & -- & $11{,}331$ & $1.34$\\
$[-0.5,-0.4)$ & $8$ & $1.00$ & $2{,}324$ & $1.93$\\
$[-0.4,-0.3)$ & $56$ & $1.73$ & $2{,}176$ & $1.94$\\
$[-0.3,-0.2)$ & $50$ & $2.82$ & $2{,}103$ & $2.42$\\
$[-0.2,-0.1)$ & $61$ & $2.93$ & $1{,}752$ & $2.33$\\
$[-0.1,0)$ & $37$ & $2.92$ & $1{,}598$ & $2.43$\\
\bottomrule
\end{tabular}
\end{table}

Mixed cycles do have one feature that $+$ cycles lack. Under $T_S$ an even number $x$ is reached
only from $(2x-1)/3$ or from $2x/3$, both smaller than $x$, so the smallest element of a $+$ cycle is
odd. In a mixed map an even $x$ can also be reached from the larger odd number $2x-1$, and some mixed
cycles have an even smallest element, for instance $(16,24,36,54,27,41,62,31)$ of
$T_{\{1,2\},[2,4]}$, in which $31$ goes down to $16$.

To sum up, the mixed maps do not seem to form a new kind. The direction at the odd numbers $x\equiv
1\pmod 4$ decides: if they go down, the map is dominated block by block by the Collatz map and, as far
as the census reaches, has only the two shifted Collatz cycles, up to switching; if they go up, the
map behaves like the maps $T_S$.

\section{Conjectures and open problems}\label{sec:open}

\begin{conjecture}\label{conj:contracting}
If $S$ is eventually periodic and contracting, then every $T_S$-orbit is eventually periodic, and
$T_S$ has finitely many cycles. For the seven maps of \cref{tab:cycles}, the cycles listed there are
all the cycles.
\end{conjecture}

For $S=\emptyset$ the first statement is the part of the Collatz conjecture that excludes divergent
orbits. For every $S$ the 2-adic integer $-1$ is a fixed point of $U_S$ (with $M(-1)=0$), and on
the residue classes that contain $-1$ every determined odd step multiplies by $3/2$. So, as for the
Collatz map, these classes have no coefficient stopping time in the sense recalled in
\cref{sec:minus}, and \cref{conj:contracting} cannot be proved from residue classes alone.

\begin{conjecture}\label{conj:expanding}
If $S\neq\N$ is expanding, then almost every $T_S$-orbit diverges. If moreover $\kappa(S)<1$, the
number of $n\le N$ with an eventually periodic orbit is $N^{1-\kappa(S)+o(1)}$.
\end{conjecture}

The condition $\kappa(S)<1$ fails only for the sets $S=\N\setminus\{c\}$. Indeed,
$\E\,2^{-X}=\frac13\sum_{m\ge 0}3^{-r_S(m)}$, and every $m\notin S$, including $m=0$, contributes $1$
to the sum. So if $S$ is expanding and misses at least two positive integers, then $\E\,2^{-X}>1$, and
$\kappa(S)<1$ by convexity; for $S=\N\setminus\{c\}$, on the other hand, $\E\,2^{-X}=1-3^{1-c}/6<1$ and
$\kappa(S)>1$. For these sets the heuristic of \cref{sec:predictions}(v) predicts that only finitely
many $n$ have an eventually periodic orbit; for $S=\N\setminus\{1\}$, for example, the cycle $(1,2)$
attracts only the numbers $1$ and $2$.

\begin{conjecture}\label{conj:basins}
If $S$ is contracting, the basin of every cycle of $T_S$ has a natural density.
\end{conjecture}

\begin{conjecture}\label{conj:branching}
$\pair{\emptyset}{\{1\}}$ is unbounded, and $\pair{\emptyset}{\{3\}}$ is bounded. More generally,
a family with $v\neq 0$ is bounded if and only if $v<0$, that is, if $\Lambda(\lambda)<1$ for some
$\lambda>0$ (\cref{sec:brw}), at least when $S^-$ and $S^+$ are eventually periodic.
\end{conjecture}

\begin{problem}
Prove that some $T_S$ with $S\neq\N$ has a divergent orbit, for instance $S=\{1,2\}$, where
$\delta=+0.377$ and almost every orbit exceeds $n^{1.18}$ (\cref{thm:density}). By
\cref{prop:rational-cycles} this requires arithmetic information beyond residue classes.
\end{problem}

\begin{problem}
Prove that every orbit of $T^-_{\{1,2\}}$ is eventually periodic. This is perhaps the simplest map in
the family that is not covered by \cref{prop:simple}: its drift $-1.21$ is about three times that of
the Collatz map, and only the steps at $x\equiv 7\pmod{16}$ go up, by a factor of about $9/8$. More
generally, for which $S$ is every $T^-_S$-orbit eventually periodic?
\end{problem}

\begin{problem}
Are $C^-_\emptyset(1)$ and $C^-_\emptyset(2)$ the only $-$ cycles, that is, does every map $T^-_S$ have
only the shifted cycles of the Collatz map? For $S=\N$ this includes the statement that the Collatz
map has no cycle other than $(1,2)$. (By \cref{prop:simple} the answer is yes for every $S$ that
satisfies its criterion, and by \cref{sec:unified-down} it is yes for cycles with smallest element
below $10^9$.) More generally, are $(1)$, $(2,3)$ and $(3,5)$ the only cycles of the unified maps
$T_{\So,\Se}$ with $1\notin\So$? Is the set of $+$ cycles infinite?
\end{problem}

\begin{problem}
Prove that $\pair{\emptyset}{\{1\}}$ is unbounded, that is, that some positive integer has a
divergent branch. By \cref{cor:nocert}, residue classes modulo $2^{10}$ do not suffice.
\end{problem}

\begin{problem}\label{prob:tao}
Extend Tao's theorem \cite{Tao22} to contracting $T_S$: is it true that for every function $f$ with
$f(n)\to\infty$ almost every $T_S$-orbit, in the sense of logarithmic density, attains a value below
$f(n)$? A similar question can be asked for a family with $v<0$: for almost every $n$, does every
branch from $n$ attain a value below $f(n)$? (The weaker statement that some branch from $n$
does so would follow from the first question, applied to $T_{S^-}$, which is contracting since
$\delta(S^-)\le v$.)
\end{problem}

\begin{problem}
By \cref{prop:nocritical} no map $T_S$ is critical. In the single-halving variant of
\cref{rem:single}, where $\rho(S)$ is replaced by $\sum_{s\in S}2^{-s}$, the critical case does occur,
for the set $S$ of positions of the binary digits $1$ of $\theta$. Describe the orbits of this map.
\end{problem}

\begin{problem}
Find other ways to transfer statements between maps. \Cref{prop:transfer} only compares families
ordered by inclusion, and \cref{rem:transfer-minus} gives a partial transfer from the maps
$T^-_{[1,b]}$ to the Collatz map. Is there, for example, a single map $T_S$ or $T^-_S$, other than
$T_\emptyset$ and $T^-_\N$, whose eventual periodicity would imply that of the Collatz map? (Among
the unified maps there are such maps, for instance $T_{[2,\infty),S}$ for every $S$, by
\cref{prop:down}(b); but they have the same block map as $T^-_\N$.)
\end{problem}

\section*{Data and code}
The program (a Python driver with a C kernel), the JSON files with all statistics, the lists of
rational witnesses of \cref{prop:rational-traps}, the script for the computation in the proof of
\cref{prop:nocritical}, 
the census programs and the lists of cycles of \cref{app:cycles} and
\cref{sec:unified}, and scripts that recompute the constants, tables and the figures are available as
ancillary files of this paper on arXiv, linked from its abstract page; the file \texttt{README.txt}
among them describes them and how to run them.

\section*{Acknowledgments}
The computations were done on Kaggle notebooks. This work was carried out with the assistance of
Claude Opus 5.5 (Anthropic), an AI system, which helped with the literature search, the derivations, the
code and the drafting of the text. The author has checked the results and takes full
responsibility for the content.

\bibliographystyle{amsplain}
\bibliography{custom}

\section*{Appendix}
\appendix
\section{Cycles across the family}\label{app:cycles}

A cycle is a finite set of numbers that a map permutes cyclically, and the same cycle can occur in
many maps of the family; here we name cycles independently of the maps. Every cycle contains odd
numbers, since a run of steps at even numbers ends at an odd number, and all odd elements of a cycle
are treated by the same rule. So a cycle is either a \emph{$+$ cycle}, whose odd elements go up to
$(3x+1)/2$, or a \emph{$-$ cycle}, whose odd elements go down to $(x+1)/2$; the former occur only for
maps $T_S$, the latter only for maps $T^-_S$. For a cycle $C$ let $U(C)$, respectively $D(C)$, be the
set of exponents $m\ge 1$ such that some element $x=2^m k$ of $C$ is mapped to $3x/2$, respectively
to $k$. These are disjoint finite sets.

\begin{proposition}\label{prop:cycle-interval}
A $+$ cycle $C$ is a cycle of $T_S$ if and only if $U(C)\subseteq S$ and $S\cap D(C)=\emptyset$, and a
$-$ cycle $C$ is a cycle of $T^-_S$ under the same condition.
\end{proposition}

\begin{proof}
$C$ is a cycle of a map exactly when the map sends every element of $C$ to its successor in $C$,
that is, when at every exponent occurring in $C$ the map follows the same rule as $C$.
\end{proof}

So the maps that contain a cycle $C$ form an interval $U(C)\subseteq S\subseteq\N\setminus D(C)$,
and this interval contains a unique map closest to the Collatz map of the same family:
$T_{U(C)}$ for a $+$ cycle, which differs from the Collatz map $T_\emptyset$ exactly at the
exponents in $U(C)$, and $T^-_{\N\setminus D(C)}$ for a $-$ cycle, which differs from the Collatz map
$T^-_\N$ exactly at the exponents in $D(C)$. By \cref{prop:threshold}(b), $T_{U(C)}$ is also the most
contracting map that contains $C$. We name a cycle by its family, by the set of exponents at which it
deviates from the Collatz rule, and by its smallest element $n$:
\begin{align*}
C^+_U(n)&:\ \text{the $+$ cycle with } U(C)=U \text{ and smallest element } n, \text{ a cycle of } T_U,\\
C^-_D(n)&:\ \text{the $-$ cycle with } D(C)=D \text{ and smallest element } n, \text{ a cycle of } T^-_{\N\setminus D}.
\end{align*}
Since $C^+_U(n)$ is the cycle of $T_U$ through $n$, and $C^-_D(n)$ the cycle of $T^-_{\N\setminus D}$
through $n$, the name determines the cycle. An empty index marks a cycle of the Collatz map:
$C^+_\emptyset(1)=(1,2)$ is a cycle of $T_\emptyset$, and $C^-_\emptyset(1)=(1)$ and
$C^-_\emptyset(2)=(2,3)$ are the Collatz cycles $(0)$ and $(1,2)$, shifted by one, as cycles of
$T^-_\N$ (\cref{prop:shift}). \Cref{tab:named} lists the cycles of all maps studied in
\cref{sec:numerics,sec:minus}.

\begin{table}[!htbp]
\centering\small
\caption{The cycles of the maps studied in \cref{sec:numerics,sec:minus}: name, the sets $U(C)$ and
$D(C)$, length in steps, number of odd steps, largest element, and the maps of
\cref{sec:numerics} (for $+$ cycles) or of \cref{sec:minus} (for $-$ cycles) that contain the
cycle; by \cref{prop:cycle-interval}, $C$ occurs exactly for the sets $S$ with $U(C)\subseteq S$ and
$S\cap D(C)=\emptyset$.}
\label{tab:named}
\begin{tabular}{lllrrrl}
\toprule
cycle & $U(C)$ & $D(C)$ & length & odd & largest & occurs for $S=$\\
\midrule
$C^+_\emptyset(1)$ & $\emptyset$ & $\{1\}$ & $2$ & $1$ & $2$ & $\emptyset,\{2\},\{3\},\even,[3,\infty),\{2,3\},[2,4]$\\
$C^+_{\{1\}}(1)$ & $\{1\}$ & $\{3\}$ & $5$ & $3$ & $8$ & $\{1\},\{1,2\}$\\
$C^+_{\{2\}}(19)$ & $\{2\}$ & $\{1\}$ & $8$ & $4$ & $66$ & $\{2\},\even,\{2,3\},[2,4]$\\
$C^+_{\{2\}}(47)$ & $\{2\}$ & $\{1,4\}$ & $24$ & $16$ & $1{,}002$ & $\{2\},\{2,3\}$\\
$C^+_{\{3\}}(3)$ & $\{3\}$ & $\{2\}$ & $4$ & $2$ & $12$ & $\{3\},[3,\infty),\odd$\\
$C^+_{\{3\}}(103)$ & $\{3\}$ & $\{1,2,4\}$ & $39$ & $27$ & $6{,}924$ & $\{3\}$\\
$C^+_{\{4\}}(15)$ & $\{4\}$ & $\{3\}$ & $6$ & $4$ & $120$ & $\even$\\
$C^+_{\{2,3\}}(5)$ & $\{2,3\}$ & $\{1,4\}$ & $18$ & $10$ & $80$ & $\{2,3\}$\\
$C^+_{[2,4]}(23)$ & $[2,4]$ & $\{1,6\}$ & $22$ & $10$ & $1{,}962$ & $[2,4]$\\
$C^+_{[3,7]}(255)$ & $[3,7]$ & $\{1,2\}$ & $126$ & $71$ & $269{,}364$ & $[3,\infty)$\\
$C^+_{[3,7]}(347)$ & $[3,7]$ & $\{1,2\}$ & $55$ & $26$ & $6{,}924$ & $[3,\infty)$\\
$C^+_{[3,9]}(271)$ & $[3,9]$ & $\{1,2\}$ & $558$ & $304$ & $7{,}409{,}556$ & $[3,\infty)$\\
\midrule
$C^-_\emptyset(1)$ & $\emptyset$ & $\emptyset$ & $1$ & $1$ & $1$ & all $S$\\
$C^-_\emptyset(2)$ & $\{1\}$ & $\emptyset$ & $2$ & $1$ & $3$ & all $S$ with $1\in S$\\
\bottomrule
\end{tabular}
\end{table}

\paragraph{A census of cycles.} To find all cycles with small elements, we search, for each $n$, the
orbits starting at $n$ in which every exponent $m\ge1$ is assigned the step up or the step down when
it first occurs, following both assignments; this is the branching family $\pair{\emptyset}{\N}$ of
\cref{sec:branching} (and its analogue with odd numbers going down), restricted to branches that
treat each exponent in the same way every time it occurs. A branch is abandoned when
it drops below $n$, when it exceeds $2^{120}$, or when it enters a loop that does not contain $n$.
This finds every cycle, of every map $T_S$ or $T^-_S$, whose smallest element is $n$ and whose largest
element is below $2^{120}$; every cycle found was checked again independently. For example, from
$n=1$ the search finds $C^+_\emptyset(1)$ (down at the exponent $1$) and $C^+_{\{1\}}(1)$ (up at $1$, down
at $3$), and it abandons the branch that goes up at $1$ and $3$ and down at $2$ when it enters the loop
$(3,5,8,12)$, which the search from $n=3$ finds as $C^+_{\{3\}}(3)$.

For $-$ cycles the bound $2^{120}$ was never reached, so the result is complete: \emph{for every
$S\subseteq\N$, the only cycles of $T^-_S$ with smallest element below $10^6$ are $C^-_\emptyset(1)$ and
$C^-_\emptyset(2)$.} (\Cref{sec:unified-down} extends this to $10^9$.) For $+$ cycles there are exactly
$53$ with smallest element below $10^5$ and
largest element below $2^{120}$; they are listed in \cref{tab:census}. Their number grows slowly with
the smallest element: $25$ of them have smallest element below $100$, $43$ below $1000$, and all of
them at most $6757$. Exactly $21$ of them occur in a contracting map, that is, have
$\rho(U(C))<\theta$. Cycles with
larger elements exist: with the bound $2^{125}$ the search finds one more cycle with smallest element
below $10^4$, $C^+_{\{1,5,9,10,11,12,13,14\}}(4033)$, of length $4371$.

\begin{figure}[!htbp]
\centering
\includegraphics[width=\textwidth]{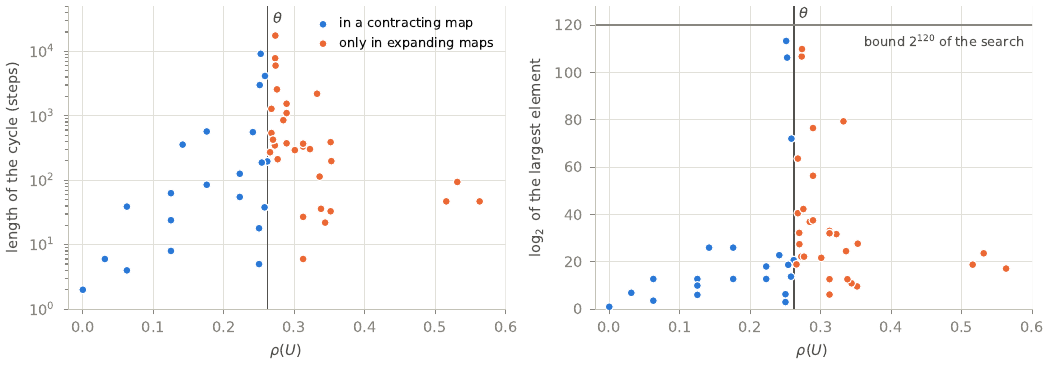}
\caption{The $53$ $+$ cycles of \cref{tab:census}: length (left) and $\log_2$ of the largest element
(right) against $\rho(U)$. Blue: cycles that occur in a contracting map ($\rho(U)<\theta$); orange:
cycles that occur only in expanding maps. The vertical line is $\theta$, the horizontal line on the
right the bound $2^{120}$ of the search.}
\label{fig:census}
\end{figure}

\Cref{fig:census} plots the $53$ cycles against $\rho(U)$. The $11$ cycles with more than $1000$ steps
and the $12$ cycles with largest element above $2^{40}$ all have $0.25<\rho(U)<0.34$, close to $\theta$
on either side, and several of them come close to the bound $2^{120}$. This fits the random model, in
which the drift is small near $\theta$ and the orbits make long excursions (\cref{tab:theory}), and it
suggests that a larger bound would reveal more cycles of this kind.

To read \cref{tab:census}, recall from \cref{prop:cycle-interval} that $C^+_U(n)$ is a cycle of
$T_S$ exactly when $U\subseteq S$ and $S\cap D=\emptyset$. For example, the cycle $C^+_{\{1,6,8,9\}}(43)$
occurs for every $S$ that contains $1,6,8,9$ and none of $2,3,4,5,7$; the most contracting of these
maps, $T_{\{1,6,8,9\}}$, has $\rho=0.2617$, just below $\theta=0.2619$.

{\footnotesize\renewcommand{\arraystretch}{0.96}
\begin{longtable}{rllrrrr}
\caption{All $+$ cycles $C^+_U(n)$ with smallest element $n<10^5$ and largest element below
$2^{120}$: $n$ (with $^*$ if the cycle occurs in a map of \cref{sec:numerics}), $U=U(C)$, $D=D(C)$,
length, number of odd steps, $\log_2$ of the largest element, and $\rho(U)$; the map $T_U$ is
contracting if and only if $\rho(U)<\theta=0.2619$.}\label{tab:census}\\
\toprule
$n$ & $U$ & $D$ & length & odd & $\log_2\max$ & $\rho(U)$\\
\midrule
\endfirsthead
\toprule
$n$ & $U$ & $D$ & length & odd & $\log_2\max$ & $\rho(U)$\\
\midrule
\endhead
\bottomrule
\endfoot
$1$$^*$ & $\emptyset$ & 1 & $2$ & $1$ & $1.0$ & $0.0000$\\
$1$$^*$ & 1 & 3 & $5$ & $3$ & $3.0$ & $0.2500$\\
$3$$^*$ & 3 & 2 & $4$ & $2$ & $3.6$ & $0.0625$\\
$5$$^*$ & 2,3 & 1,4 & $18$ & $10$ & $6.3$ & $0.2500$\\
$7$ & 1,4,5,6,7 & 2,3 & $33$ & $18$ & $9.6$ & $0.3516$\\
$9$ & 1,4,5 & 3 & $6$ & $2$ & $6.2$ & $0.3125$\\
$11$ & 1,3 & 2,5 & $27$ & $18$ & $12.7$ & $0.3125$\\
$15$$^*$ & 4 & 3 & $6$ & $4$ & $6.9$ & $0.0312$\\
$19$$^*$ & 2 & 1 & $8$ & $4$ & $6.0$ & $0.1250$\\
$21$ & 3,4,5,7 & 1,2,6 & $85$ & $50$ & $12.8$ & $0.1758$\\
$23$$^*$ & 2,3,4 & 1,6 & $22$ & $10$ & $10.9$ & $0.3438$\\
$31$ & 3,4 & 1,2,5 & $63$ & $42$ & $12.8$ & $0.1250$\\
$31$ & 1,4,5,6,7 & 2,3,9 & $390$ & $250$ & $27.7$ & $0.3516$\\
$41$ & 1,2,5,6 & 3,4,7 & $94$ & $50$ & $23.6$ & $0.5312$\\
$43$ & 1,6,8,9 & 2,3,4,5,7 & $196$ & $131$ & $20.7$ & $0.2617$\\
$43$ & 1,4,5,6,7,9 & 2,3,8 & $198$ & $119$ & $27.7$ & $0.3525$\\
$47$$^*$ & 2 & 1,4 & $24$ & $16$ & $10.0$ & $0.1250$\\
$47$ & 1,4,5 & 2,3,6 & $333$ & $212$ & $33.1$ & $0.3125$\\
$51$ & 2,3,6 & 1,5 & $38$ & $20$ & $13.7$ & $0.2578$\\
$59$ & 1,3,5,7,11 & 2,4,6,8,9,10,13 & $2197$ & $1399$ & $79.4$ & $0.3323$\\
$59$ & 2,3,6,7,8,10,11 & 1,4,5,9 & $347$ & $184$ & $22.2$ & $0.2725$\\
$61$ & 1,2,5 & 3,4,7 & $47$ & $29$ & $18.8$ & $0.5156$\\
$67$ & 2,3,6,7 & 1,4,5,8 & $272$ & $146$ & $18.9$ & $0.2656$\\
$67$ & 2,3,5,7,8,9 & 1,4,6 & $212$ & $113$ & $22.2$ & $0.2764$\\
$81$ & 1,3,6,7,8,9 & 2,5 & $36$ & $20$ & $12.7$ & $0.3379$\\
$103$$^*$ & 3 & 1,2,4 & $39$ & $27$ & $12.8$ & $0.0625$\\
$109$ & 1,2,4,5,10 & 3,9 & $47$ & $27$ & $17.2$ & $0.5630$\\
$113$ & 1,5,8 & 2,3,4,7 & $540$ & $357$ & $40.6$ & $0.2676$\\
$113$ & 1,5,9,10 & 2,3,4,6,7,8 & $1281$ & $838$ & $63.6$ & $0.2676$\\
$127$ & 2,3,7 & 1,4,5,6 & $188$ & $100$ & $18.7$ & $0.2539$\\
$215$ & 2,3,6,7,8,9 & 1,4,5,10,13 & $2568$ & $1442$ & $42.3$ & $0.2754$\\
$235$ & 1,4,9,10,11,12 & 2,3,5,6,7,8 & $854$ & $548$ & $36.9$ & $0.2844$\\
$255$$^*$ & 3,4,5,6,7 & 1,2 & $126$ & $71$ & $18.0$ & $0.2227$\\
$271$$^*$ & 3,4,5,6,7,8,9 & 1,2 & $558$ & $304$ & $22.8$ & $0.2412$\\
$281$ & 1,9 & 2,3,4,5,6,7,8,10 & $2993$ & $1987$ & $113.3$ & $0.2510$\\
$347$$^*$ & 3,4,5,6,7 & 1,2 & $55$ & $26$ & $12.8$ & $0.2227$\\
$433$ & 1,6,7,8,10,11,12,16 & 2,3,4,5,9,15 & $7804$ & $5117$ & $106.8$ & $0.2728$\\
$447$ & 1,6,7,8,10,11,12,13 & 2,3,4,5,9,16 & $17538$ & $11536$ & $106.8$ & $0.2731$\\
$463$ & 1,4,5,6 & 2,3,10 & $114$ & $72$ & $24.5$ & $0.3359$\\
$703$ & 1,3,6,8 & 2,4,5,7,9 & $304$ & $195$ & $31.7$ & $0.3223$\\
$871$ & 2,3,8,11,13,14 & 1,4,5,6,7,9,10,12 & $9143$ & $5227$ & $106.3$ & $0.2523$\\
$875$ & 1,4,6 & 2,3,5,7,8,10,16 & $1102$ & $714$ & $56.4$ & $0.2891$\\
$953$ & 1,4,6 & 2,3,5,7,8,9 & $1538$ & $991$ & $76.5$ & $0.2891$\\
$1211$ & 1,4,5,11,12 & 2,3,6,7,10 & $369$ & $240$ & $32.1$ & $0.3130$\\
$1351$ & 2,5,6,7,8 & 1,3,4,11 & $571$ & $345$ & $26.0$ & $0.1758$\\
$1711$ & 2,5,10,11 & 1,3,4,9 & $357$ & $220$ & $26.0$ & $0.1416$\\
$1817$ & 1,4,5 & 2,3,6,7,12 & $369$ & $241$ & $32.1$ & $0.3125$\\
$1819$ & 2,3,5,6,7,8 & 1,4 & $293$ & $155$ & $21.7$ & $0.3008$\\
$2559$ & 2,3,5,7 & 1,4,6,9 & $422$ & $235$ & $32.2$ & $0.2695$\\
$4033$ & 1,4,6 & 2,3,5,9 & $374$ & $239$ & $37.5$ & $0.2891$\\
$4339$ & 2,3,7,8,10,12 & 1,4,5,6,9,11 & $4144$ & $2413$ & $72.1$ & $0.2584$\\
$6503$ & 2,3,5,7,11 & 1,4,6,8,10 & $426$ & $241$ & $27.5$ & $0.2698$\\
$6757$ & 2,3,5,7,8,17,18 & 1,4,6,9,10,16 & $6003$ & $3375$ & $109.9$ & $0.2734$\\
\end{longtable}
}

\end{document}